\documentclass[11pt]{amsart}
\usepackage{amsmath}
  \usepackage{url,graphicx}
  \usepackage{amssymb, color, pstricks-add, manfnt}
  \usepackage{amsmath, amsthm,  amsfonts}
  \usepackage[plainpages=false]{hyperref}

  \theoremstyle{plain}
  \newtheorem{theorem}{Theorem}[section]
  \newtheorem{lemma}{Lemma}[section]
  \newtheorem{remark}{Remark}[section]

  \numberwithin{equation}{section}
  \numberwithin{figure}{section}

\renewcommand{\baselinestretch}{1.00}
\begin{document}

\title{On the Dirichlet problem for a class of augmented Hessian equations}

\author{Feida Jiang}
\address{Centre for Mathematics and Its Applications, The Australian National University,
              Canberra ACT 0200, Australia} \email{jfd2001@163.com}

\author{Neil S. Trudinger}
\address{Centre for Mathematics and Its Applications, The Australian National University,
              Canberra ACT 0200, Australia}
\email{Neil.Trudinger@anu.edu.au}

\author{Xiao-Ping Yang}
\address{School of Science, Nanjing University of Science and Technology, Nanjing 210094, P.R.China}
\email{yangxp@mail.njust.edu.cn}

\thanks{This research was a part of the first author's Ph.D thesis at Nanjing University of Science and Technology in April 2013. It was supported by National Natural Science Foundation of China(No.11071119) and the Australian Research Council.}


\date{\today}


\maketitle

\abstract {In this paper, we consider the Dirichlet problem for a new class of
augmented Hessian equations. Under sharp assumptions that the matrix
function in the augmented Hessian is regular and there exists a smooth
subsolution, we establish global second order derivative estimates
for the solutions to the Dirichlet problem in bounded domains.
The results extend the corresponding results in the previous paper \cite{JTY2013}
from the Monge-Amp\`ere type equations to the more general Hessian type equations.}
\endabstract


\baselineskip=12.8pt
\parskip=3pt
\renewcommand{\baselinestretch}{1.38}

\section{Introduction}\label{Section 1}
\vskip10pt

In this paper, we study a class of augmented Hessian equations with the
following form
    \begin{equation}\label{1.1}
        S_k[D^2u-A(x,Du)]=B(x,Du), \ \ {\rm in} \ \Omega,
    \end{equation}
associated with the Dirichlet boundary condition
    \begin{equation}\label{1.2}
        u=\varphi,\ \ {\rm on} \ \partial \Omega,
    \end{equation}
where $\Omega$ is a bounded domain in $n$ dimensional Euclidean
space $\mathbb{R}^n$, $A: \Omega\times\mathbb{R}^n\rightarrow
\mathbb{R}^{n}\times \mathbb{R}^{n}$ is a symmetric matrix function,
$B:\Omega\times\mathbb{R}^n\rightarrow\mathbb{R}$ is a scalar
function, $\varphi$ is a smooth function on $\partial \Omega$,
$Du$ denotes the gradient vector of $u$, $D^2u$ denotes
the Hessian matrix of the second derivatives of $u$, and $S_{k}$ is
a $k$-Hessian operator defined by
    \begin{equation}
        S_k[W]:=S_k(\lambda),
    \end{equation}
with $\lambda$ denoting the eigenvalues $\lambda_1, \cdots,
\lambda_n$ of the $n\times n$ symmetric matrix $W$, $S_k(\lambda)$
denoting the $k$-th order elementary symmetric function given by
    \begin{equation}
    S_k(\lambda)=\sum_{1\leq i_1<\cdots<i_k\leq n}\lambda_{i_1}\cdots \lambda_{i_k}, \ \ k\leq n \ {\rm is\  an \
    integer}.
    \end{equation}
In the equation \eqref{1.1}, the symmetric matrix $W$ under
consideration is the augmented Hessian matrix $\{D^2u-A(x,Du)\}$,
which is given by functions of the Hessian minus a lower order
symmetric matrix function. As usual, we shall use $(x,p)$ to define
the points in $\Omega\times \mathbb{R}^n$. We adopt the terminology
from \cite{Tru2006} and call the matrix function $A$ regular if $A$ is co-dimension one
convex with respect to $p$, in the sense that
    \begin{equation}\label{A3W}
        A_{ij,kl}(x,p)\xi_i\xi_j\eta_k\eta_l\geq 0,
    \end{equation}
for all $(x,p)\in \bar \Omega\times\mathbb{R}^n$,
$\xi,\eta\in\mathbb{R}^n$, $\xi\perp\eta$, where
$A_{ij,kl}=D^2_{p_kp_l}A_{ij}$. As explained in \cite{Tru2006},
 condition \eqref{A3W} is the natural condition for
regularity of the solution $u$.

Equations of the form \eqref{1.1} have attracted much research
interest and have many applications. If $A\equiv 0$, the equation
\eqref{1.1} reduces to the standard Hessian equation, and clearly
satisfies the regular condition \eqref{A3W}. In this case, it
becomes the standard Monge-Amp\`ere equation when $k=n$. For these
well known standard Hessian equations and Monge-Amp\`ere equations,
classical solvability of the corresponding Dirichlet problems was
studied extensively, see \cite{CNSI84}, \cite{CNSIII},
\cite{CW2001}, \cite{Guanbo1998}, \cite{Krylov1993},
\cite{Liyanyan90}, \cite{Tru95Acta}, \cite{Wang1994} and
\cite{Wang2009} etc. If $A$ depends only on $x$, the corresponding
equation \eqref{1.1} has applications in Riemannian geometry, see
\cite{Liyanyan90} and \cite{Urbas2002}. Recently, Guan studied the
Dirichlet problem for a class of Hessian type equations with $A$
only depending on $x$ on Riemannian manifolds under some very general
structure conditions in \cite{Guan2012a}.

For general $A$ related to $x$ and $p$, when $k=n$, the equation
\eqref{1.1} is the Monge-Amp\`ere type equation and is an important
model equation in optimal transportation, geometric optics,
isometric embedding and etc. For example, in optimal transportation, the matrix function $A(x,p)=D^2_{x}c(x,Y(x,p))$ and the
scalar function $B(x,p)=|\det D^2_{xy}c|\frac{f(x)}{g(Y(x,p))}$,
where $c(x,y)$ is the cost function, $Y$ is the optimal map, $f$ and
$g$ are densities of the original and target measures. The existence
of smooth solutions to the Dirichlet problem for the optimal
transportation equation in small balls under the strict version
of \eqref{A3W}, denoted A3, was used in \cite{MTW2005}.
 For more detailed definitions in optimal
transportation, one can refer to \cite{MTW2005}. Classical
solvability for Dirichlet problem of the Monge-Amp\`ere type
equations in bounded domains was studied under sharp assumptions
that $A$ is regular as well as the existence of subsolutions in
\cite{JTY2013}. While for general $k$, the augmented Hessian
equation \eqref{1.1} is related to the so-called $k$-Yamabe problem
in conformal geometry. For instance, in the conformal geometry case,
the matrix function $A$ is given by
$A(x,Du)=-\frac{1}{2}|p|^2I+p\otimes p$, see \cite{Tru2006}. One can
check that such $A$ also satisfies the regular condition
\eqref{A3W}. The Dirichlet boundary problem for a class of Hessian
type equations related to conformal deformations of metrics on
Riemannian manifolds with boundary was studied in \cite{Guanbo2007}.
Under the existence of a smooth subsolution, Guan \cite{Guanbo2007}
derived various {\it a priori} estimates and proved the classical
solvability for the Dirichlet problem.

In the current paper, we deal with the augmented Hessian equations
not only coming from conformal geometry but for the general
equation \eqref{1.1} together with the Dirichlet boundary data
\eqref{1.2}. Under the assumption that there
exists a subsolution, we show that the regular condition \eqref{A3W} is
sufficient for the second order derivative estimates of the Dirichlet problem for the
general Hessian type equation \eqref{1.1} in smooth
bounded domains without any geometric restrictions. Our main issue
is to deal with the dependence in lower order terms for both $A$ and
$B$, which was not an issue in \cite{Guan2012a}.

Before stating our theorems, we shall present some definitions and well known properties of the function
$f=(S_k)^{\frac{1}{k}}$ (see for example \cite{Guanbo2007,SUW2004,Urbas2001}). We define the positive cone $\Gamma_k^+$ in
$\mathbb{R}^n$,
    \begin{equation}
        \Gamma_k^+:=\left\{\lambda\in \mathbb{R}^n |\ S_j(\lambda)>0, \ \forall \ 1\leq j\leq
        k\right\},
    \end{equation}
and its closure,
    \begin{equation}
        \Gamma_k=\Gamma_k^+ \cup \partial \Gamma_k^+.
    \end{equation}
Then, we have
    \begin{equation}\label{1.8}
        (S_k)^{\frac{1}{k}}(\lambda)>0,\ \lambda \in \Gamma_k^+,\ \ {\rm and}\ \ (S_k)^{\frac{1}{k}}(\lambda)=0,\ \lambda\in \partial\Gamma_k^+,
    \end{equation}
    \begin{equation}\label{1.9}
        (S_k)^{\frac{1}{k}} {\ \rm is\ a \ concave \ function \ in} \ \Gamma_k^+,
    \end{equation}
    \begin{equation}\label{1.10}
        f_i:=\frac{\partial (S_k)^{\frac{1}{k}}}{\partial \lambda_i}>0,\
{\rm in} \ \Gamma_k^+,\ 1\leq i\leq n,
    \end{equation}
and for every $C>0$ and every compact set $E\subset \Gamma_k^+$,
there exists $R=R(E,C)>0$ such that
    \begin{equation}\label{1.6}
        (S_k)^{\frac{1}{k}}(\lambda_1,\cdots,\lambda_{n-1},\lambda_n+R)\geq
        C, \ \forall \  \lambda=(\lambda_1,\cdots,\lambda_n)\in E.
    \end{equation}
If $\lambda \in \Gamma_k^+$ with $\lambda_1\geq \cdots \geq
\lambda_n$, the concavity of $f$ leads to $f_1\leq \cdots \leq f_n$,
see (ii) of Lemma 2.1 in \cite{Urbas2001}. If we define $\Gamma_{k;\mu_1,\mu_2}=\{\lambda\in \Gamma_k: \ \mu_1\leq
f(\lambda)\leq \mu_2\}$ for any $0<\mu_1\leq \mu_2$, then
    \begin{equation}\label{1.12}
        \sum_{i=1}^{n}f_i\geq \sigma_0,\ {\rm on}\
        \Gamma_{k;\mu_1,\mu_2},
    \end{equation}
and, for $k>1$,
    \begin{equation}\label{1.14}
        \lim_{|\lambda|\rightarrow \infty}\sum_i f_i=+\infty,\ \
        {\rm on} \
        \Gamma_{k;\mu_1,\mu_2},
    \end{equation}
where $\sigma_0$ is a constant depending on $\mu_1$ and $\mu_2$.

The properties \eqref{1.8} and \eqref{1.9} imply the degenerate
ellipticity condition $f_i\geq 0$ in $\Gamma_k$ for $i=1,\cdots,n$.
With the definition of the cone $\Gamma_k$, we can define the admissible
solutions corresponding to the equation \eqref{1.1}. A solution
$u\in C^2(\Omega)$ is called $k$-$A$-admissible if
    \begin{equation}\label{1.18}
        D^2u-A(x,Du)\in \Gamma_k,\ \ {\rm in}\ \Omega.
    \end{equation}
For simplicity, we denote the $k$-$A$-admissible solutions as the
admissible solutions. Under the assumption $B(x,Du)>0$, the
admissible condition \eqref{1.18} ensures that the equation
\eqref{1.1} is elliptic with respect to a solution $u\in
C^2(\Omega)$.
Denoting the right hand side function by
$\tilde B(x,p)=(B(x,p))^{\frac{1}{k}}$, we assume that $\tilde B$
is convex with respect to the gradient variables $p$, that is
    \begin{equation}\label{convexity on B}
        D^2_{p_kp_l}\tilde B(x,p)\xi_k\xi_l \geq 0,
    \end{equation}
for all $(x,p)\in \Omega\times \mathbb{R}^n$, $\xi\in R^n$.

We now state our main theorems as follows:

\begin{theorem}\label{Th1}
Let $u\in C^4(\Omega)\cap C^2(\bar\Omega)$ be an admissible solution
of Dirichlet problem \eqref{1.1}-\eqref{1.2}, where $A\in
C^2(\bar\Omega\times\mathbb{R}^n)$ is regular, $B\in
C^2(\bar\Omega\times\mathbb{R}^n)$ satisfies \eqref{convexity on B} and $B>0$ in
 $\bar\Omega\times\mathbb{R}^n$. Suppose also there
exists an admissible function $\underline u\in C^2(\bar\Omega)$  satisfying
    \begin{equation}\label{Gamma+}
        D^2\underline u-A(x,D\underline u)\in \Gamma_k^+,\ \  { for \
        all }\  x\in \bar\Omega.
    \end{equation}
Then we have the estimate
    \begin{equation}\label{interior bound}
        \sup\limits_\Omega |D^2u| \leq C(1+\sup\limits_{\partial
        \Omega} |D^2u|),
    \end{equation}
where the constant $C$ depends on $k, n, A, B, \Omega, \underline u$
and $\sup\limits_\Omega{(|u|+|Du|)}$.
\end{theorem}

We remark that the function $\underline u$ in this theorem is only assumed to satisfy \eqref{Gamma+}.  Note that Theorem \ref{Th1} certainly holds if $\underline u$ is an admissible subsolution satisfying
    \begin{equation}\label{sub}
        S_k[D^2\underline u-A(x,D\underline u)]\geq B(x,D\underline u), \ \ {\rm in} \ \Omega.
    \end{equation}
The convexity condition \eqref{convexity on B} on $\tilde B$ can be removed in the Monge-Amp\`ere case when $k=n$, see \cite{JTY2013}. When $k=1$, the conclusion \eqref{interior bound} follows from
 the classical Schauder theory, \cite{GT2001}.

From Theorem \ref{Th1}, we can infer a global second derivative bound
for solutions of the Dirichlet problem \eqref{1.1}-\eqref{1.2} from boundary estimates,
which we can  further derive  if $\underline u$ is
an admissible subsolution satisfying \eqref{sub} with the given boundary trace.
Therefore, we have the following theorem.
\begin{theorem}\label{Th2}
In addition to the assumptions in Theorem \ref{Th1}, suppose the admissible function $\underline u$ is a
subsolution satisfying \eqref{sub} in $\Omega$ and $\underline u = \varphi$ on $\partial\Omega$ with
$\varphi \in C^4(\partial \Omega), \partial\Omega\in C^4$. Then any
admissible solution $u\in C^4(\Omega)\cap C^2(\bar\Omega)$ of the
Dirichlet problem \eqref{1.1}-\eqref{1.2} satisfies the global
{\it a priori} estimate
    \begin{equation}\label{second bound}
        \sup\limits_\Omega |D^2u|\leq C,
    \end{equation}
where the constant $C$ depends on $k, n, A, B, \Omega, \varphi,
\underline u$ and $\sup\limits_\Omega{(|u|+|Du|)}$.
\end{theorem}

In these theorems, the {\it a priori} estimates for second
derivatives are derived under hypotheses including the regular
condition for $A$ and the existence of a subsolution. We remark that these
assumptions are sharp for the second order derivative estimates to the
Dirichlet problem \eqref{1.1}-\eqref{1.2} in general bounded domains.
Moreover, if $A$ and $B$ depend also on $u$, we have the more
general augmented Hessian equations,
    \begin{equation}\label{general}
        S_k[D^2u-A(x,u,Du)]=B(x,u,Du), \ \ {\rm in} \ \Omega.
    \end{equation}
We shall also discuss the corresponding theorems in Section \ref{Section 5} under
some mild additional structure conditions for both $A$ and $B$ with
respect to $u$.

We focus on the second order derivative estimates in the current paper.
These second derivative estimates together with the solution bounds
and the gradient estimates will yield the regularity and the classical
existence results of the Dirichlet problem \eqref{1.1}-\eqref{1.2}
for the augmented Hessian equations. Once we obtain the derivative estimates up to second order,
the Evans-Krylov theorems then yields $C^{2,\alpha}$ bounds of the solutions.
Then we obtain the existence and uniqueness of the classical solution of the Dirichlet problem by the method of continuity.
It would be interesting to derive the lower order estimates especially the gradient
estimates under some additional structure conditions analogous to the natural conditions
of Ladyzhenskaya and Ural'tseva for quasilinear elliptic equations\cite{LU1968,Tru83,GT2001}. This will be taken up in a sequel.

The paper is organized as follows: In Section \ref{Section 2}, we introduce some
preliminary lemmas. The first lemma is to construct a global barrier
function for the linearized operator of $F$, which is fundamental in
the second derivative estimates. In Section \ref{Section 3}, we prove the global
second derivative estimates for solutions to equation \eqref{1.1},
which reduce the second derivative bound to the boundary. In Section
\ref{Section 4}, we first show that the regular condition of the matrix $A$ is
preserved when we translate and rotate the coordinates. Then we
obtain the second derivative bound on the boundary for regular $A$. In Section
\ref{Section 5}, we discuss more general equations \eqref{general}, where both
$A$ and $B$ depend also on $u$. More structure conditions on $A$ and
$B$ will be explained for this general case.

\section{Preliminaries}\label{Section 2}

In this section, we introduce some notation and present some
preliminary results needed in later sections.

We denote the augmented Hessian matrix by $W$, that is
    \begin{equation}\label{2.1}
        W=\{w_{ij}\}=\{u_{ij}-A_{ij}(x,Du)\}.
    \end{equation}
Let
    \begin{equation}\label{2.2}
        F[u]=:F(w_{ij})=(S_{k})^{\frac{1}{k}}(w_{ij})=(S_k)^{\frac{1}{k}}[u_{ij}-A_{ij}(x,Du)],
    \end{equation}
it is known that $F$ is a concave operator with respect to $w_{ij}$ for admissible $u$.
Introducing the linearized operator of $F$:
    \begin{equation}\label{2.3}
        L=F^{ij}[D_{ij}-D_{p_k}A_{ij}(x,Du)D_k],
    \end{equation}
where $F^{ij}=\frac{\partial F}{\partial w_{ij}}$, it follows that $\{F^{ij}\}$ is positive definite \cite{CNSIII}. We also set
    \begin{equation}\label{2.4}
        \mathcal {L}=L-\tilde B_{p_i}D_i=F^{ij}[D_{ij}-D_{p_k}A_{ij}(x,Du)D_k]-\tilde B_{p_i}D_i.
    \end{equation}

We shall first prove the following fundamental lemma, which is
a key barrier construction needed in the second derivative
estimates. Although its proof is similar to Lemma 2.1 in
\cite{JTY2013}, for completeness and for the convenience in later
sections, we still present the detailed proof. Also we take the opportunity to make a small correction
to our previous proof, (in connection with the choice of $x_0$).
\begin{lemma}\label{Lemma 2.1}
Let $u\in C^2(\bar \Omega)$ be an admissible solution of equation
\eqref{1.1}, $\underline u\in C^2(\bar \Omega)$ be an admissible strict
subsolution of equation \eqref{1.1} satisfying
    \begin{equation}\label{subsolution}
        S_k[D^2\underline u-A(x,D\underline u)]\geq B(x,D\underline
        u)+\delta_0,
    \end{equation}
for some positive constant $\delta_0$. Assume the matrix function
$A$ is regular satisfying \eqref{A3W}, $A_{ij}(x,p)\in
C^2(\bar\Omega\times\mathbb{R}^n)$, $i,j=1,\cdots,n$, $B(x,p)\in
C^2(\bar\Omega\times\mathbb{R}^n)$. Then
    \begin{equation}\label{2.5}
        \mathcal {L} \left({e^{K(\underline u-u)}}\right) \geq
        \epsilon_1\sum\limits_{i} F^{ii}-C,
    \end{equation}
holds in $\Omega$ for positive constants $K$, $\epsilon_1$ and $C$, which depend on $k, A, B, \Omega, \sup\limits_\Omega|Du|$
and $\sup\limits_\Omega|D \underline u|$.
\end{lemma}
\noindent \textbf{Proof}. Since $\underline u$ is a strict
subsolution of \eqref{1.1}, for any $x_0\in \Omega$, the perturbation function $\underline u_\epsilon=\underline u -\frac{\epsilon}{2}|x-x_0|^2$ is still a strict
subsolution, for sufficiently small $\epsilon > 0$, and satisfies
    \begin{displaymath}
        F[\underline u_\epsilon]=(S_{k})^{\frac{1}{k}} [D^2\underline
        u_\epsilon-A(x,D\underline u_\epsilon)]\geq (B(x,D\underline
        u_\epsilon)+\tau)^{\frac{1}{k}},
    \end{displaymath}
for some positive constant $\tau$.

Let $v=\underline u-u$, $v_\epsilon=\underline u_\epsilon-u$. By a direct
calculation, we have
    \begin{equation}\label{2.6}
        \begin{array}{lll}
            Lv&=& L(v_\epsilon)+L(\frac{\epsilon}{2}|x-x_0|^2)\\
            &=& \epsilon F^{ii}-\epsilon
            F^{ij}D_{p_k}A_{ij}(x,Du)(x-x_0)_k +F^{ij}\{D_{ij}(\underline u_\epsilon -u)-[A_{ij}(x,D\underline
            u_\epsilon)-A_{ij}(x,Du)]\}\\
            & & +F^{ij}\{A_{ij}(x,D\underline u_\epsilon)-A_{ij}(x,Du)-D_{p_k}A_{ij}(x,Du)D_kv_\epsilon\}.
        \end{array}
    \end{equation}
Since $F$ is concave with respect to $w_{ij}$, we have
    \begin{displaymath}
        F[\underline u_\epsilon]-F[u]\leq F^{ij}\{D_{ij}(\underline
        u_\epsilon -u)-[A_{ij}(x,D\underline u_\epsilon)-A_{ij}(x,Du)]\}.
    \end{displaymath}
Using the Taylor expansion, we have for some $\theta\in (0,1)$,
    \begin{equation}\label{2.7}
        \begin{array}{ll}
            &A_{ij}(x,D\underline
            u_\epsilon)-A_{ij}(x,Du)-D_{p_k}A_{ij}(x,Du)D_kv_\epsilon\\
            =&D_{p_k}A_{ij}(x,\hat{p})D_kv_\epsilon-D_{p_k}A_{ij}(x,Du)D_kv_\epsilon\\
            =&\frac{\theta}{2}A_{ij,kl}(x,\bar p)D_kv_\epsilon
            D_lv_\epsilon,\\
        \end{array}
    \end{equation}
where $\hat{p}=(1-\theta)Du+\theta D\underline u_\epsilon$, $\bar
p=(1-\bar \theta)Du+\bar \theta D\underline u_\epsilon$ and $\bar
\theta\in (0, \theta)$.
Next we choose a finite family of balls $B_\rho(x^i)$, with centres at $x^i$, $ i = 1\cdots N$, covering $\bar\Omega$
and with fixed radii $\rho < 1/2n\max\limits_{\bar \Omega} |D_{p_k}A_{ij}(x,Du)|$. Then we select $x_0 = x^i$ for some $i$ such that $x\in B_\rho(x^i)$.
Accordingly, we have for a fixed positive $\epsilon$,
    \begin{equation}\label{fixed epsilon}
        \begin{array}{lll}
            Lv & \geq &  \epsilon F^{ii}-\epsilon
            F^{ij}D_{p_k}A_{ij}(x,Du)(x-x_0)_k+F[\underline
            u_\epsilon]-F[u]+\frac{\theta}{2}F^{ij}A_{ij,kl}(x,\bar p)D_kv D_lv\\
            & \geq &  \epsilon F^{ii}-\frac{\epsilon}{2n}|F^{ij}|+\frac{\theta}{2}F^{ij}A_{ij,kl}(x,\bar
            p)D_kv D_lv+ (B(x,D\underline u)+\tau)^{\frac{1}{k}}-(B(x,Du))^{\frac{1}{k}}\\
            & \geq &  \epsilon F^{ii}-\frac{\epsilon}{2n}|F^{ij}|+\frac{\theta}{2}F^{ij}A_{ij,kl}(x,\bar
            p)D_kv D_lv -C_1,
        \end{array}
    \end{equation}
holds for $x\in B_\rho(x^i)$, where $C_1$ is a constant depending on $B$, $Du$, and $D\underline u$.
We see that \eqref{fixed epsilon} holds in all balls $B_\rho(x^i)$, $i=1\cdots N$,
with a fixed positive constant $\epsilon$. Then by the finite covering,
\eqref{fixed epsilon} holds in $\Omega$ with a uniform positive constant $\epsilon$.

Let $\phi=e^{Kv}$ with positive constant $K$ to be determined, we
have
    \begin{displaymath}
        \begin{array}{lll}
            L\phi & = & Ke^{Kv}Lv + K^2e^{Kv}F^{ij}D_ivD_jv\\
            & \geq & \displaystyle Ke^{Kv}\left\{\epsilon
            F^{ii}-\frac{\epsilon}{2n}|F^{ij}|+\frac{\theta}{2}F^{ij}A_{ij,kl}(x,\bar p)D_kv D_lv
            -C_1+KF^{ij}D_ivD_jv\right\}.
        \end{array}
    \end{displaymath}
Without loss of generality, assume that $Dv=(D_1v, 0, \cdots,
0)$, we get
    \begin{displaymath}
        \begin{array}{lll}
            L\phi & \geq & Ke^{Kv}\left\{\epsilon
            F^{ii}-\frac{\epsilon}{2n}|F^{ij}|+\frac{\theta}{2}F^{ij}A_{ij,11}(x,\bar p)(D_1v)^2
            +KF^{11}(D_1v)^2-C_1\right\}\\
            & \geq & Ke^{Kv}\left\{\epsilon
            F^{ii}-\frac{\epsilon}{2n}|F^{ij}|+\frac{\theta}{2}\sum\limits_{i \hspace{1mm} or \hspace{1mm} j
            = 1}F^{ij}A_{ij,11}(x,\bar p)(D_1v)^2
            +KF^{11}(D_1v)^2-C_1\right\},\\
        \end{array}
    \end{displaymath}
here we use the fact that $A$ is regular in the second inequality.

Since the matrix $\{F^{ij}\}$ is positive definite, any $2\times 2$
diagonal minor has positive determinant. By the Cauchy's inequality, we have
    \begin{displaymath}
        |F^{ij}|\leq \sqrt{F^{ii}F^{jj}}\leq \frac{1}{2}(F^{ii}+F^{jj}), \ \
         {\rm (which\  leads \ to} \sum\limits_{i,j}|F^{ij}|\leq n\sum\limits_i F^{ii} {\rm)},
    \end{displaymath}
and
    \begin{displaymath}
        |F^{1i}|\leq \sqrt{F^{11}F^{ii}}\leq \eta F^{ii}+\frac{1}{4\eta}F^{11},
    \end{displaymath}
for any positive constant $\eta$.

Thus, we have
    \begin{displaymath}
        \begin{array}{lll}
            L\phi & \geq & \displaystyle Ke^{Kv}\left\{\frac{\epsilon}{2}
            F^{ii}-\theta \eta F^{ii}|A_{1i,11}(x,\bar
            p)|(D_1v)^2-\frac{\theta}{4\eta}F^{11}|A_{1i,11}(x,\bar p)|(D_1v)^2
            +KF^{11}(D_1v)^2-C_1\right\}.\\
        \end{array}
    \end{displaymath}
Choosing $\eta$ small such that $\eta \leq \frac{\epsilon}{4\theta
\max{\{|A_{1i,11}(x,\bar p)|(D_1v)^2\}}}$ and $K$ large such that
$K\geq \frac{\theta \max|A_{1i,11}(x,\bar p)|}{4\eta}$, we obtain
    \begin{displaymath}
        \begin{array}{lll}
            L\phi & \geq & \displaystyle Ke^{Kv}\{
            \frac{\epsilon}{4}F^{ii}-C_1\}.\\
        \end{array}
    \end{displaymath}
Thus, we have
    \begin{displaymath}
        \mathcal {L}\phi=L\phi-\tilde B_{p_i}D_i\phi\geq Ke^{Kv}\{
        \frac{\epsilon}{4}F^{ii}-C_1\}-\tilde B_{p_i}D_i\phi.
    \end{displaymath}
If we choose
$\epsilon_1=\min\limits_{\bar\Omega}\{\frac{\epsilon}{4}K e^{Kv}\}$
and $C=\max\limits_{\bar\Omega}\{C_1Ke^{Kv}+\tilde
B_{p_i}D_i\phi\}$, the conclusion of this lemma is proved.

\qed

\begin{remark}
The global barrier construction in this lemma is fundamental to
derive the second order derivative estimates for solutions to
\eqref{1.1}-\eqref{1.2}. In Section \ref{Section 3}, we shall use this barrier
function to reduce the global estimates of second order derivatives
to boundary estimates when the matrix function $A$ is regular. If
$A$ is strictly regular satisfying \eqref{A3} in the next section,
we do not need such a barrier function, as further
discussed in Remark \ref{remark 3.1}. In Section \ref{Section 4}, this barrier function will be
modified a bit to fit the barrier argument so that we can get the
second order derivative bounds on the boundary for both the mixed
tangential-normal and the double normal directions.
\end{remark}

\begin{remark}
In previous papers \cite{Tru2006,TruWang2009}, another global
barrier condition called $A$-boundedness condition is assumed,
namely that a domain $\Omega$ is $A$-bounded with respect to $u$, if
there exists a function $\varphi\in C^2(\bar \Omega)$ satisfying
    \begin{equation}\label{Abound}
        [D_{ij}\varphi - D_{p_k}A_{ij}(x,Du)D_k
        \varphi]\xi_i\xi_j\geq \delta_0 |\xi|^2,
    \end{equation}
for some $\delta_0>0$ and for all $x\in \Omega$, $\xi\in
\mathbb{R}^n$. When the diameter of $\Omega$ is sufficiently small,
the function $\varphi=|x|^2$ satisfies condition \eqref{Abound} for
bounded $Du$. Also, condition \eqref{Abound} is trivial in the
standard Monge-Amp\`ere case as seen by still taking
$\varphi(x)=|x|^2$. In the optimal transportation case, there are also
various examples showing that condition \eqref{Abound} is satisfied
by regular cost functions, see \cite{LiuTru2010}. From the condition
\eqref{Abound}, we immediately have the following inequality
    \begin{equation}
        \mathcal{L}\varphi = F^{ij}[D_{ij}\varphi - D_{p_k}A_{ij}(x,Du)D_k
        \varphi]-\tilde B_{p_i}D_i\varphi \geq \delta_0 F^{ii} - C,
    \end{equation}
which has a similar form of the inequality \eqref{2.5} in Lemma
\ref{Lemma 2.1}. In this sense, the $A$-boundedness condition
\eqref{Abound} can provide us with an alternative global barrier
function $\varphi$.
\end{remark}

\begin{remark}\label{remark 2.3}
By adding the perturbation function $a e^{bx_1}$ for small positive
constant $a$ and large positive constant $b$, a non-strict classical
subsolution for an elliptic partial differential equation can be
made strict using the linearized operator and the mean value
theorem, (see \cite{GT2001}, Chapter 3). This can also be done near
the boundary, preserving the boundary condition by adding the
perturbation $a e^{bd(x)}$, where $a$ is small positive constant,
$b$ is a large positive constant, and $d(x)=dist(x,\partial \Omega)$
is the distance function. Hence we need only assume the existence of
a non-strict subsolution in Lemma 2.1; the inequality (\ref{2.5})
will still hold for the corresponding strict subsolution. Thus, the
second order {\it{apriori}} estimates will also hold under the
existence of a non-strict subsolution and we only need to assume a
non-strict subsolution in the hypotheses for our theorems.
\end{remark}

\begin{remark}\label{remark 2.4}
We observe from the proof of Lemma \ref{Lemma 2.1} that the function $\underline u$ does not need to be a subsolution.
If $\underline u$ is only an admissible function satisfying $D^2\underline u-A(x,D\underline u)\in \Gamma_k^+$ in $\bar \Omega$,
then we have $F[\underline u]=(S_k)^{\frac{1}{k}}(D^2\underline u-A(x,D\underline u))\geq \delta^{\frac{1}{k}}>0$ for some $\delta>0$. The corresponding perturbation function $\underline u_\epsilon=\underline u -\frac{\epsilon}{2}|x-x_0|^2$ still satisfies $F[\underline u_\epsilon]=(S_k)^{\frac{1}{k}}(D^2\underline u_\epsilon-A(x,D\underline u_\epsilon))\geq \tau$ for some positive constant $\tau$, and the conclusion \eqref{2.5} of Lemma \ref{Lemma 2.1} still holds. Hence, the existence of a non-strict subsolution in Remark \ref{remark 2.3} can be further relaxed by the existence of an admissible function satisfying $D^2\underline u-A(x,D\underline u)\in \Gamma_k^+$ in $\bar\Omega$. Furthermore in the optimal transportation case such functions are constructed in \cite{TruWang2009} whence, in particular, the subsolution condition can be removed altogether from Theorems 1.1 and 2.1  in \cite{JTY2013}, for equations arising from optimal transportation; (see \cite{JT2014}). This also facilitates a more direct proof of Theorem 1.1 in \cite{Tru2013}.
\end{remark}

As usual we denote the second partial derivatives of $F$ with
respect to $w_{ij}$ by $F^{ij,kl}$, that is
$F^{ij,kl}=\frac{\partial^2 F}{\partial w_{ij}\partial w_{kl}}$. We also
need the following lemma of Andrews \cite{Andrews94, Ger96}, for our
global estimates in the next section.

\begin{lemma}\label{Lemma 2.2}
For any $n\times n$ symmetric matrix $W=\{w_{ij}\}$, one has that
    \begin{equation}
        F^{ij,kl}w_{ij}w_{kl}=\sum_{i,j}\frac{\partial^2f}{\partial \lambda_i\partial
        \lambda_j}w_{ii}w_{jj}+\sum_{i\neq
        j}\frac{f_i-f_j}{\lambda_i-\lambda_j}w^2_{ij}.
    \end{equation}
The second term on the right hand side is non-positive if $f$ is
concave, and is interpreted as a limit if $\lambda_i = \lambda_j$.
\end{lemma}

\section{Global second derivative estimates}\label{Section 3}

In this section, by using the preparatory lemmas in Section 2, we
give the proof of Theorem 1.1, which shows that  the global bounds for second
derivative estimates of the equation \eqref{1.1} are reduced to
their boundary estimates. Our arguments mimic those in
\cite{MTW2005} and \cite{TruWang2009}, and are modifications of the
arguments presented in Section 17.6 of \cite{GT2001}. In the following proof, the function $\underline u$ can be regarded
as an admissible function satisfying \eqref{Gamma+}, as explained in Remark \ref{remark 2.4}.

\vspace{3mm}

\noindent \textbf{Proof of Theorem 1.1}. Let $v$ be an auxiliary
function given by
    \begin{equation}
        v(x,\xi):=\log w_{\xi\xi}+ \eta(\frac{1}{2}|Du|^2)+b\phi,
    \end{equation}
where $w_{\xi\xi}=w_{ij}\xi_i\xi_j=(u_{ij}-A_{ij})\xi_i\xi_j$ with a
vector $\xi\in \mathbb{R}^n$, $\phi=e^{K(\underline u-u)}$ is the
barrier function as in Lemma \ref{Lemma 2.1} with $\underline u$ satisfying \eqref{Gamma+}, $\eta$ is a function
to be determined and $b$ is a positive constant to be determined.

By differentiating the following equation in the $\xi$ direction,
    \begin{equation}\label{F[u]}
        F[u]=:F(u_{ij}-A_{ij}(x,Du))=(S_k)^{\frac{1}{k}}[u_{ij}-A_{ij}(x,Du)]=\tilde B(x,Du),
    \end{equation}
we have,
    \begin{equation}\label{differentiate once}
        \begin{array}{lll}
            F^{ij}D_\xi w_{ij}& = & F^{ij}\left[D_{ij}u_\xi-D_\xi
            A_{ij}-(D_{p_k}A_{ij})D_ku_\xi\right]\\
            & = & D_\xi\tilde{B}+(D_{p_k}\tilde{B})D_ku_\xi.
        \end{array}
    \end{equation}
By a further differentiation, we obtain
    \begin{equation}\label{differentiate twice}
        \begin{array}{lll}
            F^{ij}D_{\xi\xi}w_{ij}& = & \displaystyle
            F^{ij}\left[D_{ij}u_{\xi\xi}-D_{\xi\xi}A_{ij}-2(D_{\xi
            p_k}A_{ij})D_ku_{\xi}-
            (D^2_{p_kp_l}A_{ij})D_ku_{\xi}D_lu_{\xi}-(D_{p_k}A_{ij})D_ku_{\xi\xi}\right]\\
            &= & \displaystyle -F^{ij,kl}D_{\xi}w_{ij}D_{\xi}w_{kl}+
            D_{\xi\xi}\tilde{B}+2(D_{\xi p_k}\tilde{B})
            D_ku_{\xi}+(D^2_{p_kp_l}\tilde{B})D_ku_{\xi}D_lu_{\xi}+(D_{p_k}\tilde{B})D_ku_{\xi\xi}.
        \end{array}
    \end{equation}

Assume that $v$ takes its maximum at an interior point $x_0 \in
\Omega$ and a unit vector $\xi_0$, without loss of generality we can
choose an orthogonal coordinate system $e_1,\cdots,e_n$ at this
point such that $e_1(x_0)=\xi_0$, $\{w_{ij}\}$ is diagonal and
$w_{11}\geq \cdots \geq w_{nn}$. We immediately have $F^{ij}$ is
diagonal and $F^{11}\leq \cdots \leq F^{nn}$. At the point $x_0$,
the function $v(x,e_1)=\log w_{11}+ \eta(\frac{1}{2}|Du|^2)+b\phi$
attains its maximum. By a direct calculation, we have at the maximum point
$x_0$,
    \begin{equation}\label{differentiate once '}
        D_iv=\frac{D_iw_{11}}{w_{11}}+\eta'D_kuD_{ik}u+bD_i\phi=0,
    \end{equation}
    \begin{equation}
        D_{ij}v=\frac
        {D_{ij}w_{11}}{w_{11}}-\frac{D_iw_{11}D_jw_{11}}{w^2_{11}}+\eta'(D_{ik}uD_{jk}u+D_kuD_{ijk}u)+\eta''D_kuD_{ik}uD_{l}uD_{jl}u+bD_{ij}\phi \leq 0,
    \end{equation}
and
    \begin{equation}\label{Lv leq 0}
        \begin{array}{rll}
            0 \geq \mathcal {L}v & = &
            F^{ij}[D_{ij}v-(D_{p_k}A_{ij})D_kv]-(D_{p_k}\tilde{B})D_kv \\
            & = & \displaystyle \frac{1}{w_{11}}\mathcal {L}w_{11}+\eta'\sum\limits_ku_k\mathcal
            {L}u_k+b\mathcal {L}\phi-\frac{1}{w^2_{11}}F^{ii}(D_iw_{11})^2+\eta'F^{ii}(D_{ki}u)^2+\eta''F^{ii}(D_kuD_{ki}u)^2.
        \end{array}
    \end{equation}

We shall estimate each term in (\ref{Lv leq 0}), we first
observe that the term $\eta'\sum\limits_ku_k\mathcal {L}u_k$ has a lower bound by using (\ref{differentiate once}), that is
    \begin{equation}
        \eta'\sum u_k\mathcal {L}u_k \geq -\eta'_0C (\sum F^{ii}+1),
    \end{equation}
where $\eta'_0=|\eta'|_{C^0}$, the positive constant $C$ depends on
$n, A, B, \Omega$ and $\sup\limits_\Omega{(|u|+|Du|)}$. Unless
otherwise specified, we shall use $C$ to denote a positive constant
with such dependance in this section.

By the barrier construction in Lemma \ref{Lemma 2.1}, we have an
estimate for $b\mathcal{L}\phi$, that is
    \begin{equation}
        b\mathcal{L}\phi= b\mathcal{L}\left(e^{K(\underline u-u)}\right)\geq
        b\epsilon_1 \sum F^{ii}-bC,
    \end{equation}
where the positive constants $\epsilon_1$ and $C$ also depend on $\underline u$.

By the twice differentiated equation (\ref{differentiate twice}), we
have
    \begin{equation}
        \begin{array}{rl}
            \mathcal{L}u_{11}=
            & -F^{ij,kl}D_1w_{ij}D_1w_{kl}+F^{ij}(D_{11}A_{ij}+2(D_{1p_k}A_{ij})u_{k1})\\
            &+F^{ij}D^2_{p_kp_l}A_{ij}u_{k1}u_{l1}+D_{11}\tilde{B}+2(D_{1p_k}\tilde{B})u_{k1}+(D^2_{p_kp_l}\tilde{B})u_{k1}u_{l1},
        \end{array}
    \end{equation}
where
    \begin{displaymath}
        F^{ij}(D_{11}A_{ij}+2(D_{1p_k}A_{ij})u_{k1})=F^{ij}(D_{11}A_{ij}+2(D_{1p_k}A_{ij})(w_{k1}+A_{k1}))\geq
        -CF^{ii}(1+w_{11}),
    \end{displaymath}
and
    \begin{displaymath}
        \begin{array}{rl}
            F^{ij}D^2_{p_kp_l}A_{ij}u_{k1}u_{l1} & =
            F^{ij}A_{ij,kl}(w_{k1}+A_{k1})(w_{l1}+A_{l1})\\
            & =
            F^{ij}(A_{ij,kl}w_{k1}w_{l1}+2A_{ij,kl}w_{k1}A_{l1}+A_{ij,kl}A_{k1}A_{l1})\\
            & =
            F^{ij}(A_{ij,11}w^2_{11}+2A_{ij,1l}w_{11}A_{l1}+A_{ij,kl}A_{k1}A_{l1})\\
            & \geq F^{11}A_{11,11}w^2_{11}-CF^{ii}(1+w_{11})\\
            & \geq -CF^{11}w^2_{11}-CF^{ii}(1+w_{11}),
        \end{array}
    \end{displaymath}
here we use the regular condition \eqref{A3W} for the matrix $A$ in the first
inequality. So we obtain that
    \begin{equation}\label{Lu11}
        \begin{array}{lll}
            \mathcal{L}u_{11} & \geq & -F^{ij,kl}D_1w_{ij}D_1w_{kl}
            +(D^2_{p_kp_l}\tilde{B})u_{k1}u_{l1}-
            C(1+w_{11}+F^{ii}+F^{ii}w_{11})-CF^{11}w^2_{11}\\
            & \geq & -F^{ij,kl}D_1w_{ij}D_1w_{kl} - C[F^{ii}(1+w_{11})]-CF^{11}w^2_{11},
        \end{array}
    \end{equation}
here we use property \eqref{1.12}: $\sum f_i\geq \sigma_0$ on
$\Gamma_{k;\mu_1,\mu_2}$ and the convexity condition \eqref{convexity on B} to derive the second inequality.

To obtain the estimate $\mathcal {L}w_{11}$, we also need to
estimate the term $\mathcal {L}A_{11}$. We evaluate this term in
detail,
    \begin{displaymath}
        \begin{array}{rll}
            \mathcal {L}A_{11} & = &
            F^{ij}[D_{ij}A_{11}+D_{ip_k}A_{11}(w_{kj}+A_{kj})+D_{jp_k}A_{11}
            (w_{ki}+A_{ki})\\
            &   &
            +D^2_{p_kp_l}A_{11}(w_{ki}+A_{ki})(w_{lj}+A_{lj})+D_{p_k}A_{11}u_{kij}-D_{p_k}A_{ij}D_kA_{11}\\
            &   &
            -D_{p_k}A_{ij}D_{p_l}A_{11}(w_{kl}+A_{kl})]-D_{p_k}\tilde{B}[D_kA_{11}+D_{p_l}A_{11}(w_{kl}+A_{kl})]\\
            & \leq &
            C+CF^{ii}(1+w_{jj})+F^{ij}D^2_{p_kp_l}A_{11}w_{ki}w_{lj}+F^{ij}D_{p_k}A_{11}u_{kij}\\
            & \leq &
            CF^{ii}(1+w_{jj})+CF^{ii}w^2_{ii}+F^{ij}D_{p_k}A_{11}u_{kij}.\\
        \end{array}
    \end{displaymath}
The above third derivative term can be estimated by using
the differentiated equation (\ref{differentiate once}), that is
    \begin{displaymath}
        \begin{array}{rll}
            F^{ij}D_{p_k}A_{11}u_{kij}& = & F^{ij}D_{p_k}A_{11}u_{ijk}=
            D_{p_k}A_{11}\{F^{ij}[A_{ijk}+D_{p_l}A_{ij}u_{lk}]+\tilde{B}_k+D_{p_l}\tilde{B}u_{lk}\}\\
            & = & D_{p_k}A_{11}\{F^{ij}[A_{ijk}+D_{p_l}A_{ij}(w_{lk}+A_{lk})]+\tilde{B}_k+D_{p_l}\tilde{B}(w_{lk}+A_{lk})\}\\
            & \leq & C[F^{ii}(1+w_{jj})+w_{ii}],
        \end{array}
    \end{displaymath}
here the property \eqref{1.12}: $\sum f_i\geq \sigma_0$ on
$\Gamma_{k; \mu_1,\mu_2}$ is used. So
we get
    \begin{equation}\label{LA11}
        \mathcal {L}A_{11} \leq C[F^{ii}(1+w_{jj})+w_{ii}]+CF^{ii}w^2_{ii}.
    \end{equation}
Combining (\ref{Lu11}) and (\ref{LA11}), we obtain the estimate for
$\mathcal {L}w_{11}$,
    \begin{equation}\label{Lw11}
        \mathcal {L}w_{11} = \mathcal {L}u_{11}-\mathcal {L}A_{11} \geq
        -F^{ij,kl}D_1w_{ij}D_1w_{kl} - C[F^{ii}(1+w_{jj})+w_{ii}]-CF^{ii}w^2_{ii}-CF^{11}w^2_{11}.
    \end{equation}

With the above estimates in hand, \eqref{Lv leq 0} becomes
    \begin{equation}\label{Lv leq 0 '}
        \begin{array}{lll}
            0\geq\mathcal {L}v \geq & \displaystyle -\frac{1}{w_{11}}F^{ij,kl}D_1w_{ij}D_1w_{kl}
            -\frac{1}{w^2_{11}}F^{ii}(D_iw_{11})^2\\
            & \displaystyle
            +\eta'F^{ii}(D_{ik}u)^2+\eta''F^{ii}(D_kuD_{ik}u)^2-\frac{C}{w_{11}}F^{ii}w^2_{ii}-CF^{11}w_{11}\\
            & \displaystyle
            -\frac{C}{w_{11}}[F^{ii}(1+w_{jj})+w_{ii}]+(b\epsilon_1-\eta'_0C)F^{ii}-C(b+\eta'_0).
        \end{array}
    \end{equation}
Firstly, we need to estimate the first two terms on the right hand side of \eqref{Lv leq 0 '}. We recall that
$w_{11}$ is the largest eigenvalue of $\{w_{ij}\}$, and set
    \begin{equation}
    I=\{i:w_{ii}\leq -\theta w_{11}\},\ \ J=\{i>1:w_{ii}>-\theta w_{11}\},
    \end{equation}
where $\theta\in (0,1)$ is a positive constant to be chosen later.
We have $1 \notin I$, $1 \notin J$, $I\cap J=\emptyset$, $\{1\}\cup
I \cup J=\{1,2,\cdots,n\}$, and $J^c={1}\cup I$, where $J^c$ is the
complementary set of $J$. By Lemma \ref{Lemma 2.2} and the concavity
of the operator $F=(S_k)^{\frac{1}{k}}$, we have
    \begin{equation}
        \begin{array}{lll}
        \displaystyle -\frac{1}{w_{11}}F^{ij,kl}D_1w_{ij}D_1w_{kl}& \geq & \displaystyle -\frac{1}{w_{11}}\sum_{i\neq
        j}\frac{F^{ii}-F^{jj}}{w_{jj}-w_{ii}}(D_1w_{ij})^2\\
        & \geq & \displaystyle -\frac{2}{w_{11}}\sum_{i\geq 2}\frac{F^{ii}-F^{11}}{w_{11}-w_{ii}}(D_1w_{i1})^2\\
        & \geq & \displaystyle -\frac{2}{w_{11}}\sum_{i\in J}\frac{F^{ii}-F^{11}}{w_{11}-w_{ii}}(D_1w_{i1})^2\\
        & \geq & \displaystyle \frac{2}{(1+\theta)w^2_{11}}\sum_{i\in J}(F^{ii}-F^{11})[D_iw_{11}+(D_iA_{11}-D_1A_{i1})]^2\\
        & \geq & \displaystyle \frac{2-2\theta}{(1+\theta)w^2_{11}}\sum_{i\in J}(F^{ii}-F^{11})[(D_iw_{11})^2-\frac{1}{\theta}(D_iA_{11}-D_1A_{i1})^2]\\
        & \geq & \displaystyle \frac{1}{w^2_{11}}\sum_{i\in J}F^{ii}(D_iw_{11})^2-\frac{C}{w^2_{11}}\sum_{i\in J}F^{ii}-\frac{F^{11}}{w^2_{11}}\sum_{i\in J}(D_iw_{11})^2.
        \end{array}
    \end{equation}
here we use the Cauchy's inequality in the last second inequality,
and we chose $\theta= 1/3$ to obtain the last inequality. Therefore,
by the relationship between the set $I$ and $J$, we have
    \begin{equation}\label{key estiamte}
        \begin{array}{ll}
                 &\displaystyle -\frac{1}{w_{11}}F^{ij,kl}D_1w_{ij}D_1w_{kl}-\frac{1}{w^2_{11}}F^{ii}(D_iw_{11})^2\\
            \geq &\displaystyle \sum_{i\in I}F^{ii}(\frac{D_iw_{11}}{w_{11}})^2-\frac{C}{w^2_{11}}\sum_{i\in J}F^{ii}-2F^{11}\sum_{i\notin I}(\frac{D_iw_{11}}{w_{11}})^2\\
            \geq &\displaystyle \sum_{i\in I}F^{ii}(\eta'D_kuD_{ik}u+bD_i\phi)^2-\frac{C}{w^2_{11}}\sum F^{ii}-2F^{11}\sum_{i\notin
            I}(\eta'D_kuD_{ik}u+bD_i\phi)^2\\
            \geq &\displaystyle (\eta')^2\sum_{i\in I}F^{ii}(D_kuD_{ik}u)^2-b^2C\sum_{i\in I}F^{ii}-\frac{C}{w^2_{11}}\sum F^{ii}-CF^{11}[(\eta')^2w^2_{11}+(\eta')^2+b^2].
        \end{array}
    \end{equation}
where the second inequality is from \eqref{differentiate once '}.

Next, if we choose the function $\eta(t)=\frac{a}{2}(1+t)^2$, we have
$\eta'(t)=a(t+1)$, $\eta''(t)=a$ and
$\eta''(t)-(\eta')^2=a-a^2(1+t)^2$. For any $t\in [0,C]$,
we choose the positive constant $a$ sufficiently small such that $\eta''(t)-(\eta')^2\geq 0$.
Therefore, we have now determined the function
    \begin{equation}\label{eta}
        \eta(\frac{1}{2}|Du|^2)=\frac{a}{2}(1+\frac{1}{2}|Du|^2)^2,
    \end{equation}
where $a$ is a small positive constant.

Consequently, by \eqref{Lv leq 0 '}, \eqref{key estiamte} and
\eqref{eta}, we have
    \begin{equation}\label{Lv leq 0 ''}
        \begin{array}{lll}
            0\geq\mathcal {L}v & \geq & \displaystyle \eta'F^{ii}(D_{ik}u)^2+[\eta''-(\eta')^2]F^{ii}(D_kuD_{ik}u)^2-\frac{C}{w_{11}}F^{ii}w^2_{ii}-CF^{11}w_{11}\\
            && \displaystyle
            -\frac{C}{w_{11}}[F^{ii}(1+w_{jj})]+(b\epsilon_1-\eta'_0C)F^{ii}-C(b+\eta'_0)\\
            && \displaystyle
            -b^2C\sum_{i\in I}F^{ii}-\frac{C}{w^2_{11}}\sum
            F^{ii}-CF^{11}[(\eta')^2w^2_{11}+(\eta')^2+b^2]\\
            &\geq & \displaystyle
            (\eta'-\frac{C}{w_{11}})F^{ii}w^2_{ii}-b^2C\sum_{i\in I}F^{ii}+(b\epsilon_1-\eta'_0C-C)F^{ii}-C(b+\eta'_0)\\
            && \displaystyle
            -CF^{11}[(\eta')^2w^2_{11}+(\eta')^2+b^2]-CF^{11}w_{11}.\\
        \end{array}
    \end{equation}
Note that we can always
suppose $w_{11}$ as large as we want, otherwise $w_{11}$ is bounded
and the proof is finished. We first choose $w_{11}$ large such that $\eta'/2\geq C/w_{11}$.
By the property \eqref{1.14}, when $k>1$, we
can choose the constant $b$ sufficiently large such that
$(b\epsilon_1-\eta'_0C-C)F^{ii}-C(b+\eta'_0)> 0$. Thus, we obtain
from \eqref{Lv leq 0 ''},
    \begin{equation}\label{eta 'Fii 1}
        \displaystyle \frac{\eta'}{2}F^{ii}w^2_{ii}\leq b^2C\sum_{i\in
        I}F^{ii}+C(\eta')^2F^{11}w^2_{11}+Cb^2F^{11}+CF^{11}w_{11}.
    \end{equation}
Since we have
    \begin{equation}\label{eta 'Fii 2}
        \begin{array}{ll}
            \displaystyle \frac{\eta'}{2}F^{ii}w^2_{ii}&\displaystyle \geq \frac{\eta'}{2}F^{11}w^2_{11} +\frac{\eta'}{2}\sum_{i\in I}F^{ii}w^2_{ii}\\
            &\displaystyle \geq \frac{\eta'}{2}F^{11}w^2_{11} +\frac{\eta'}{18}w^2_{11}\sum_{i\in
            I}F^{ii}.
        \end{array}
    \end{equation}
By choosing $a>0$ sufficiently small and $w_{11}$ sufficiently large, from \eqref{eta 'Fii 1} and \eqref{eta 'Fii 2}, we obtain the estimate, at $x_0$
    \begin{equation}\label{3.22}
        \frac{a}{4} F^{11}w^2_{11}\leq Cb^2F^{11}+CF^{11}w_{11},
    \end{equation}
that is
    \begin{equation}\label{3.23}
        \frac{a}{4} w^2_{11}\leq Cb^2+Cw_{11},
    \end{equation}
which directly leads to
    \begin{equation}
        w_{11}\leq C.
    \end{equation}
This implies the conclusion \eqref{interior bound} for $2 \leq k \leq n$. While $k=1$, since we have $F^{ii}=1$ for $i=1\cdots n$, then conclusion \eqref{interior bound} can be easily derived from \eqref{Lv leq 0 ''}. Note that in this case the conclusion (1.22) also follows from the classical Schauder theory in \cite{GT2001}. Therefore, we complete the proof of Theorem \ref{Th1}.

\qed

\begin{remark}\label{remark 3.1}
We call the matrix $A$ strictly regular if
    \begin{equation}\label{A3}
        A_{ij,kl}(x,p)\xi_i\xi_j\eta_k\eta_l >
        |\xi|^2|\eta|^2,
    \end{equation}
 for all
$(x,p)\in\bar\Omega\times\mathbb{R}^n$, $\xi,\eta\in\mathbb{R}^n$,
$\xi\perp\eta$. The global second derivative estimate
\eqref{interior bound} for the solution $u$ will become much simpler
if $A$ is strictly regular. In this case, we only need to consider
the simpler auxiliary function $v:=w_{\xi\xi}$ and maximize it over
$\bar \Omega$, since the strict regular condition directly implies
    \begin{equation}
        F^{ij}D_{p_kp_l}A_{ij}D_{k\xi}uD_{l\xi}u\geq C_1
        F^{ii}|D^2u|^2,
    \end{equation}
for some $C_1>0$. By differentiating the equation \eqref{F[u]}
twice, we then have
    \begin{equation}
        \begin{array}{ll}
        0\geq \mathcal{L}v & \geq F^{ii}[C_1|D^2u|^2-C_2(|D^2u|+1)]+(D^2_{p_kp_l}\tilde{B})u_{k1}u_{l1}\\
                           & \geq F^{ii}[C'_1|D^2u|^2-C_2(|D^2u|+1)],
        \end{array}
    \end{equation}
holds at the maximum point $x_0$, for positive constants $C_1$,
$C'_1$ and $C_2$, here the property \eqref{1.14} is used to obtain
the second inequality. Therefore, one easily obtains the global
second derivative estimate \eqref{interior bound}. Thus, to obtain
the global second derivative estimate, the barrier construction in
Lemma \ref{Lemma 2.1} is not required
 when the matrix function $A$ is strictly regular. Moreover in this case, the convexity
condition \eqref{convexity on B} on $B$ is also not needed
 and we obtain more general interior estimates;
(see  \cite{Tru2006}, Theorem 2.1).
\end{remark}

\begin{remark}
Under the assumption that the matrix $A$ is regular satisfying \eqref{A3W},
we need the convexity condition \eqref{convexity on B} on the right hand side $\tilde B$
to obtain the global second derivative estimates for both the Hessian operator $(S_k)^{\frac{1}{k}}$ here and the Hessian quotient operator ${(S_k/S_l)}^{\frac{1}{k-l}}$ ($l<k\leq n$), see \cite{vonNessi2010} for the special case ${(S_n/S_l)}^{\frac{1}{n-l}}$ ($l<n$); while in the Monge-Amp\`ere case, the special structure of the determinant function
allows us to obtain the second derivative estimate even if $\tilde B$ does not satisfy the convexity condition, see \cite{TruWang2009,LiuTru2010,JTY2013}.
If an operator $f$ satisfies a stronger property, namely $\sum\limits_i f_i/ |\lambda|\rightarrow \infty$ as $|\lambda|\rightarrow \infty$, we also do not need to impose the
convexity condition on the right hand side function $\tilde B$, see \cite{Urbas1995,Urbas2001} for references.
\end{remark}

\section{Boundary estimates for second derivatives}\label{Section 4}

In this section, we shall establish the second derivative estimate $|D^2u|\leq C$ on
the boundary $\partial \Omega$ and finish the proof of Theorem
\ref{Th2}. First, we note from Remark \ref{remark 2.3} that a non-strict
subsolution can be made strict near the boundary. Hence, we can
assume the subsolution $\underline u$ is strict provided we restrict
to a neighbourhood of $\partial \Omega$. Accordingly, Lemma
\ref{Lemma 2.1} can be retained in a neighbourhood of $\partial
\Omega$, which will suffice for the boundary estimates. Note that in this section,
$\underline u$ denotes the subsolution rather than merely an admissible function satisfying \eqref{Gamma+}.

Next, we need to check the invariance properties of the equation
\eqref{1.1} under translation and rotation of coordinates. For a
fixed point $x_0\in \bar\Omega$, let $\breve x = x - x_0$, the
equation \eqref{1.1} becomes
\[S_k[D^2_{\breve x}u - A(\breve x + x_0,D_{\breve x}u)]= B(\breve x+x_0,D_{\breve x}u).\]
Let $\breve A (\breve x,D_{\breve x} u)= A(\breve x+x_0,D_{\breve
x}u)$ and $\breve B (\breve x,D_{\breve x} u)= B(\breve
x+x_0,D_{\breve x}u)$, we have
\[S_k[D^2_{\breve x}\breve u - \breve A(\breve x,D_{\breve x} u)] = \breve B(\breve x,D_{\breve x} u).\]
Accordingly we see that both the form of the equation and the regular
condition \eqref{A3W} are invariant under translation of
coordinates. So, we may suppose any $x_0\in \bar \Omega$ as the
origin if necessary.

Let $\hat x = R x$, where $R$ is a rotation matrix. Since the
rotation matrix $R$ is orthogonal satisfying $det R =1$, the equation
\eqref{1.1} becomes
\[S_k\left\{R^T[D^2_{\hat x}u - RA(R^{-1}\hat x,R^{-1}D_{\hat x}u)R^{-1}]R\right\}= B(R^{-1}\hat x,R^{-1}D_{\hat x}u).\]
Observing that $S_k$ is invariant under rotation of coordinates, we
have
\[S_k[D^2_{\hat x}u - \hat A(\hat x, D_{\hat x}u)]= \hat B(\hat x, D_{\hat x}u),\]
with
    \begin{equation*}
        \left\{
            \begin{array}{ll}
                \hat A(\hat x, D_{\hat x}u)=RA(R^{-1}\hat x,R^{-1}D_{\hat
                x}u)R^{-1},\\
                \hat B(\hat x, D_{\hat x}u)=B(R^{-1}\hat x,R^{-1}D_{\hat x}u).
            \end{array}
        \right.
    \end{equation*}
Thus, we obtain the invariance of both the form of
equation \eqref{1.1} and the regular condition \eqref{A3W} under
rotations of coordinates.


Consequently, for any given boundary point $x_0\in
\partial\Omega$, by a translation and a rotation of the coordinates, we may take $x_0$ as the
origin and take the positive $x_n$ axis to be the inner normal of
$\partial \Omega$ at the origin. Near the origin, $\partial \Omega$
can be represented as a graph
    \begin{equation}
        x_n=\rho(x'),
    \end{equation}
such that $D\rho(0)=0$, where $x'=(x_1,\cdots,x_{n-1})$. Since
$u-\underline u =0$ on $\partial \Omega$, we then have by
differentiation twice,
    \begin{equation}\label{Dalphabeta}
        D_{\alpha\beta}(u-\underline u)(0)=-D_n(u-\underline
        u)(0)\rho_{\alpha\beta}(0), \ \ \alpha,\beta=1,\cdots,n-1,
    \end{equation}
which leads to the double tangential derivative estimate
$|D_{\alpha\beta}u(0)|\leq C$, $\alpha,\beta=1,\cdots,n-1$.

We then estimate the mixed tangential-normal derivatives $D_{\alpha
n}u(0)$, $\alpha=1,\cdots,n-1$. To obtain this estimation, we shall
apply the standard barrier argument by modifying the barrier
function introduced in Lemma 2.1. Rewrite equation \eqref{1.1} in
the following form,
    \begin{equation}\label{3.2}
        F[u]=:F(u_{ij}-A_{ij}(x,Du))=(S_k)^{\frac{1}{k}}[u_{ij}-A_{ij}(x,Du)]=\tilde
        B(x,Du).
    \end{equation}
By differentiating \eqref{3.2} with respect to $x_k$, we have
    \begin{equation}\label{3.3}
        F^{ij}(u_{k i j}-D_kA_{ij}(x,Du)-D_{p_l}A_{ij}(x,Du)u_{k l})=\tilde
        B_k(x,Du)+\tilde B_{p_l}(x,Du)u_{k l}, \ \ k=1, \cdots, n,
    \end{equation}
which leads to
    \begin{equation}\label{LDku}
        \mathcal {L}D_k u =\tilde B_k (x,Du) + F^{ij}D_kA_{ij}(x,Du),\ \ k=1,
        \cdots, n,
    \end{equation}
where $\mathcal {L}$ is defined by (\ref{2.4}). For fixed $\alpha<n$,
we consider the following operator
    \begin{equation}
        T=\partial_\alpha + \sum_{\beta<n}
        \rho_{\alpha\beta}(0)(x_{\beta}\partial_n-x_n\partial_{\beta}).
    \end{equation}
By calculation, we have
    \begin{equation}\label{LTu}
        \begin{array}{lll}
            \mathcal{L}Tu & = & \displaystyle \mathcal{L}[D_\alpha u+\sum_{\beta <n}\rho_{\alpha\beta}(0)(x_\beta D_nu-x_n D_\beta
            u)]\\
            & = & \displaystyle \mathcal{L}D_\alpha u+\sum_{\beta <n}\rho_{\alpha\beta}(0)(x_\beta \mathcal{L}D_nu-x_n \mathcal{L}D_\beta
            u)\\
            & & \displaystyle + \sum_{\beta
            <n}\rho_{\alpha\beta}(0)\left [2(F^{\beta j}u_{nj}-F^{nj}u_{\beta j})-F^{ij}(D_{p_\beta}A_{ij}u_n-D_{p_n}A_{ij}u_\beta)-(\tilde B_{p_{\beta}}u_n-\tilde
            B_{p_{n}}u_\beta)\right ].
        \end{array}
    \end{equation}
For fixed $\beta<n$, we observe that
    \begin{equation}\label{Fbeta}
        F^{\beta j}u_{nj}
        =F^{\beta j}A_{nj}\ \ {\rm and}\ \  F^{nj}u_{\beta j}=F^{nj}A_{\beta j}.
    \end{equation}
Combining \eqref{LTu}, \eqref{Fbeta}, and the differentiated
equation \eqref{LDku} for $k=1,\cdots,n-1$, we derive
    \begin{equation}
        |\mathcal{L}T(u-\underline u)|\leq C(1+\sum_i F^{ii}).
    \end{equation}
We also observe that, on $\partial \Omega$ near the origin,
    \begin{equation}
        |T(u-\underline u)|\leq C|x|^2.
    \end{equation}

Next, we shall modify the barrier function constructed in Lemma 2.1
to get through barrier argument near the boundary. Let $d=d(x)$ be
the distance function from $\partial \Omega$, we may take $\delta>0$
small enough so that $d$ is a smooth function in
$\Omega_\delta=\Omega\cap B_{\delta}(0)$. The key ingredient is the
following lemma:
\begin{lemma}\label{Lemma 3.1}
For any $M>0$, there exist constants $K$, $N$ sufficiently large and
$\mu$, $\delta$ sufficiently small, such that the function
\[\psi=1-\exp \left\{K\left[(\underline u-u)-\mu d+Nd^2\right]\right\},\]
satisfies
    \begin{equation}\label{3.6}
        \mathcal {L}\psi\leq -\frac{\epsilon_1}{2}\sum\limits_{i} F^{ii}-M, \ \   {\rm in}\ \  \Omega_\delta,\ \ \  {\rm and}\ \ \
        \psi\geq 0, \ \  {\rm on}\ \  \partial \Omega_\delta,
    \end{equation}
for some positive constant $\epsilon_1$.
\end{lemma}

\noindent \textbf{Proof}. We shall follow the proof of Lemma \ref{Lemma 2.1} and make some
necessary changes. By calculation, we have
    \begin{equation}\label{3.7}
        |Ld|\leq C\sum_i F^{ii},
    \end{equation}
and
\[Ld^2=2dLd+2F^{ij}d_id_j,\]
where $L$ is the linearized operator defined by \eqref{2.3}.

For any $x_0\in \Omega_\delta$, the perturbation function
$\underline u_\epsilon=\underline u -\frac{\epsilon}{2}|x-x_0|^2$ is
still a strict subsolution of the equation \eqref{1.1}. Let
$v=(\underline u-u)-\mu d+Nd^2$, by a calculation as in \eqref{2.6}, we
have
    \begin{equation}\label{3.8}
        \begin{array}{ll}
            &Lv=L[(\underline u-u)-\mu d+Nd^2]\\
            =& \epsilon F^{ii}-\epsilon
            F^{ij}D_{p_k}A_{ij}(x,Du)(x-x_0)_k +F^{ij}[D_{ij}\underline u_\epsilon -A_{ij}(x,D\underline
            u_\epsilon)+2N d_id_j]\\
            &-F^{ij}[D_{ij}u-A_{ij}(x,Du)]+F^{ij}[A_{ij}(x,D\underline
            u_\epsilon)-A_{ij}(x,Du)-D_{p_k}A_{ij}(x,Du)D_kv_\epsilon]\\
            &+(2Nd-\mu)Ld.
        \end{array}
    \end{equation}
By the concavity of $F$, we have
    \begin{equation}\label{3.9}
        \begin{array}{ll}
            &F^{ij}[D_{ij}\underline u_\epsilon -A_{ij}(x,D\underline
            u_\epsilon)+2N
            d_id_j]-F^{ij}[D_{ij}u-A_{ij}(x,Du)]\\
            \geq & F(D_{ij}\underline u_\epsilon -A_{ij}(x,D\underline
            u_\epsilon)+2N
            d_id_j)-F(D_{ij}u-A_{ij}(x,Du)).
        \end{array}
    \end{equation}
Since $d_n(0)=1$, $d_\beta(0)=0$ for all $\beta<n$, and $d$ is a
smooth function in $\Omega_\delta$, we can choose $\delta$
sufficiently small such that $d_n\geq 1/2$ and $d_\beta\leq 1/2$ in
$\Omega_\delta$. By the property \eqref{1.6}, we observe that the
term $F(D_{ij}\underline u_\epsilon -A_{ij}(x,D\underline
u_\epsilon)+2Nd_id_j)$ can be made as large as we want by choosing
$N$ large enough. If $N$ is chosen large enough such that
    \begin{equation}\label{3.10}
        F(D_{ij}\underline u_\epsilon -A_{ij}(x,D\underline
u_\epsilon)+2Nd_id_j)\geq \max\limits_{\bar \Omega}\tilde B(x,Du)+
\tilde M,
    \end{equation}
where $\tilde M$ is a large positive constant to be determined.
We fix the radius $\delta<1/2n\max |D_{p_k}A_{ij}(x,Du)|$. By
\eqref{2.7}, \eqref{3.7} and \eqref{3.10}, for $x\in \Omega_\delta$, we have from \eqref{3.8},
    \begin{equation}
        \begin{array}{ll}
            Lv & \geq \epsilon F^{ii}-\epsilon F^{ij}D_{p_k}A_{ij}(x,Du)(x-x_0)_k+\frac{\theta}{2}F^{ij}A_{ij,kl}(x,\bar
            p)D_kv D_lv - C(\mu+2N\delta) F^{ii} +\tilde M\\
               & \geq \frac{\epsilon}{2} F^{ii}+\frac{\theta}{2}F^{ij}A_{ij,kl}(x,\bar
            p)D_kv D_lv - C(\mu+2N\delta) F^{ii} +\tilde M.
        \end{array}
    \end{equation}
Let $\phi=e^{Kv}$ with positive constant $K$ to be determined.
Following the same way to estimate $L\phi$ in Lemma \ref{Lemma 2.1},
(by using the regularity condition \eqref{A3W} of the matrix
function $A$ and the Cauchy's inequality), we have
    \begin{equation}
        L\phi \geq Ke^{Kv}\{\frac{\epsilon}{4}F^{ii}-C(\mu +
        2N\delta)F^{ii}+\tilde M\},
    \end{equation}
holds for sufficiently large positive constant $K$, here $K$ is
fixed now.

Thus, by choosing $\mu={\epsilon}/{16C}$ and
$\delta={\epsilon}/{32NC}$ small enough, we have
    \begin{equation}
        \mathcal {L}\phi=L\phi-\tilde B_{p_i}D_i\phi\geq
        Ke^{Kv}\{\frac{\epsilon}{8}F^{ii}+\tilde M\}-\tilde B_{p_i}D_i\phi.
    \end{equation}
We choose $\tilde M = \max\limits_{\bar\Omega}\{\frac{M+\tilde
B_{p_i}D_i\phi}{Ke^{Kv}}\}$, and set
$\epsilon_1=\min\limits_{\bar\Omega}\{\frac{\epsilon}{4}K e^{Kv}\}$,
(note that $Ke^{Kv}$ is small when $K$ is a large constant,
$\epsilon_1$ is the same constant as in Lemma \ref{Lemma 2.1}), then
    \begin{equation}
        \mathcal {L}\phi\geq \frac{\epsilon_1}{2} \sum_i F^{ii} + M.
    \end{equation}
Since the constant $\tilde M$ is now fixed, the large enough
constant $N$ can be fixed now to guarantee the inequality
\eqref{3.10}. Therefore, we have
    \begin{equation}\label{3.15}
        \mathcal {L}\psi=\mathcal {L}(1-\phi)\leq -\frac{\epsilon_1}{2}
        \sum_i F^{ii} - M, \ \ {\rm in} \ \Omega_\delta.
    \end{equation}

On $\partial \Omega \cap \Omega_\delta$, we have $\psi=0$. Since
$\mu$ and $\delta$ are given by $\mu={\epsilon}/{16C}$ and
$\delta={\epsilon}/{32NC}$ respectively, we get
    \begin{equation}
        v= (\underline u-u)-\mu d+Nd^2\leq (-\mu +N\delta)d= -\frac{\epsilon d}{32C}\leq
        0,\hspace{3mm}{\rm on}\hspace{3mm}\Omega\cap\partial\Omega_\delta.
    \end{equation}
Consequently, we have $\phi=e^{Kv} \leq 1$ on
$\Omega\cap\partial\Omega_\delta$, furthermore,
    \begin{equation}
       \psi=1-\phi\geq 0,\hspace{3mm}{\rm on}\hspace{3mm}\Omega\cap\partial\Omega_\delta.
    \end{equation}
Therefore, we get $\psi\geq 0$ on $\partial \Omega_\delta$. Together
with the inequality \eqref{3.15}, the conclusion \eqref{3.6} of
Lemma \ref{Lemma 3.1} is proved.

\qed

With the barrier function $\psi$ in hand, we shall employ a new
barrier function
    \begin{equation}\label{3.18}
        \tilde \psi =a \psi + b|x|^2,
    \end{equation}
with positive $a$ and $b$ to be determined. We can choose $a \gg b
\gg 1$, then
    \begin{equation}
        \mathcal{L}\tilde \psi =a\mathcal{L}\psi + b\mathcal{L}|x|^2\leq
        -\frac{a\epsilon_1}{2}(1+\sum_i F^{ii}),
    \end{equation}
and
    \begin{equation}
        \tilde \psi\geq b|x|^2,\ \ {\rm on} \ \partial \Omega_\delta.
    \end{equation}
Therefore, for $a\gg b\gg 1$ and $\delta \ll 1$, there holds
    \begin{equation}
        \left\{
            \begin{array}{rllll}
                |\mathcal{L}T(u-\underline u)|+\mathcal{L}\tilde \psi & \leq & 0,           & {\rm in} & \Omega_\delta,\\
                |T(u-\underline u)|              & \leq & \tilde \psi, & {\rm on} & \partial \Omega_\delta.
            \end{array}
        \right.
    \end{equation}
By the maximum principle, we obtain the mixed tangential-normal
derivative estimate
    \begin{equation}
        |D_{\alpha n}u(0)|\leq C,\ \  \alpha=1,\cdots,n-1.
    \end{equation}

Finally, we will adapt the technique in \cite{Tru95Acta} to estimate the
double normal derivative $D_{nn}u$ on the boundary. Note that the
regularity of the matrix function $A$ is also critical for this
estimate. We observe that equation \eqref{1.1} can be rewritten as
    \begin{equation}\label{3.25}
        (D_{nn}u-A_{nn}(x,Du))S_{k-1}\left\{[D^2u-A(x,Du)]'\right\}+R=B(x,Du),\
    \end{equation}
where
$[D^2u-A(x,Du)]'=\{D_{\alpha\beta}-A_{\alpha\beta}(x,Du)\}_{1\leq
\alpha,\beta\leq n-1}$, $R$ denotes the remaining terms which do not involve
$D_{nn}u-A_{nn}(x,Du)$. From the above double tangential
derivative and mixed tangential-normal derivative bounds, the term
$R$ is bounded. Since an upper bound
    \begin{equation}
        D_{nn}u(x)\leq C,
    \end{equation}
is equivalent to an upper bound
    \begin{equation}
        D_{nn}u(x)-A_{nn}(x,Du)\leq C,
    \end{equation}
in order to obtain an upper bound $D_{nn}u(x)\leq C$, by
\eqref{3.25} we only need to get a positive lower bound for
$S_{k-1}\{[D^2u-A(x,Du)]'\}$.

For any boundary point $x\in \partial \Omega$, let
$\xi^{(1)},\cdots,\xi^{(n-1)}$ be an orthogonal vector field on
$\partial \Omega$. Writing
    \begin{equation*}
        \nabla_\alpha u=\xi^{(\alpha)}_mD_mu,\ \ \nabla_{\alpha\beta}u=\xi^{(\alpha)}_m\xi^{(\beta)}_lD_{ml}u,\ \
        \mathfrak{C}_{\alpha\beta}=\xi^{(\alpha)}_m\xi^{(\beta)}_lD_m\gamma_l,\
        \ 1\leq \alpha,\beta\leq n-1,
    \end{equation*}
    \begin{equation*}
        \nabla u= (\nabla_1u,\cdots,\nabla_{n-1}u),\  \
        \mathcal{A}_{\alpha\beta}(x,\nabla
        u,-D_{\gamma}u)=\xi^{(\alpha)}_m\xi^{(\beta)}_lA_{ml}(x,\nabla
        u,-D_{\gamma}u),\ \ 1\leq \alpha,\beta\leq n-1,
    \end{equation*}
and
    \begin{equation*}
        \nabla^2u=\{\nabla_{\alpha\beta}u\}_{1\leq \alpha,\beta\leq n-1},\ \ \mathfrak{C}=\{\mathfrak{C}_{\alpha\beta}\}_{1\leq \alpha,\beta\leq n-1},
    \end{equation*}
where $\gamma$ is the unit outer normal of $\partial \Omega$ at $x$.
From the boundary condition $u=\varphi$ on $\partial \Omega$, we
have
    \begin{equation}
        \nabla^2(u-\varphi)=D_\gamma(u-\varphi)\mathfrak{C},
    \end{equation}
which agree with \eqref{Dalphabeta}. Thus, we have, on the boundary
$\partial \Omega$,
    \begin{equation}
        \nabla_{\alpha\beta}u-\mathcal{A}_{\alpha\beta}(x,\nabla u,-D_\gamma u)
        =D_\gamma(u-\varphi)\mathfrak{C}_{\alpha\beta}+\nabla_{\alpha\beta}\varphi-\mathcal{A}_{\alpha\beta}(x,\nabla \varphi,-D_\gamma u), \ \ \alpha,\beta=1,\cdots,n-1,
    \end{equation}
Denote
    \begin{equation}
        \begin{array}{ll}
            G[u](x):&=G(\nabla_{\alpha\beta}u-\mathcal{A}_{\alpha\beta}(x,\nabla u,-D_\gamma u))\\
            &=G(D_\gamma(u-\varphi)\mathfrak{C}_{\alpha\beta}+\nabla_{\alpha\beta}\varphi
            -\mathcal{A}_{\alpha\beta}(x,\nabla \varphi,-D_\gamma u))\\
            &=\left\{S_{k-1}[D_\gamma(u-\varphi)\mathfrak{C}_{\alpha\beta}+\nabla_{\alpha\beta}\varphi
            -\mathcal{A}_{\alpha\beta}(x,\nabla \varphi,-D_\gamma
            u)]\right\}^{\frac{1}{k-1}}\ \ {\rm on} \ \partial
            \Omega,
        \end{array}
    \end{equation}
and
    \begin{equation}
        G^{\alpha\beta}=\frac{\partial G}{\partial r_{\alpha\beta}},
    \end{equation}
where
$r_{\alpha\beta}=D_\gamma(u-\varphi)\mathfrak{C}_{\alpha\beta}+\nabla_{\alpha\beta}\varphi
            -\mathcal{A}_{\alpha\beta}(x,\nabla \varphi,-D_\gamma u)$.

Assume that $\inf\limits_{x\in \partial \Omega}G$ is attained at
$x_0\in \partial \Omega$, we have for all $x\in \partial \Omega$,
    \begin{equation}\label{3.29}
    \begin{array}{ll}
        & G(D_\gamma(u-\varphi)(x)\mathfrak{C}_{\alpha\beta}(x)+\nabla_{\alpha\beta}\varphi(x)
            -\mathcal{A}_{\alpha\beta}(x,\nabla \varphi(x),-D_\gamma u(x))) \\
        \geq &
        G(D_\gamma(u-\varphi)(x_0)\mathfrak{C}_{\alpha\beta}(x_0)+\nabla_{\alpha\beta}\varphi(x_0)
            -\mathcal{A}_{\alpha\beta}(x_0,\nabla \varphi(x_0),-D_\gamma u(x_0))).
    \end{array}
    \end{equation}
By the concavity of $G$, we have
    \begin{equation}\label{Gconcavity}
    \begin{array}{ll}
        & G_{x_0}^{\alpha\beta}[D_\gamma(u-\varphi)(x)\mathfrak{C}_{\alpha\beta}(x)+\nabla_{\alpha\beta}\varphi(x)
            -\mathcal{A}_{\alpha\beta}(x,\nabla \varphi(x),-D_\gamma u(x))] \\
        \geq &
        G_{x_0}^{\alpha\beta}[D_\gamma(u-\varphi)(x_0)\mathfrak{C}_{\alpha\beta}(x_0)+\nabla_{\alpha\beta}\varphi(x_0)
            -\mathcal{A}_{\alpha\beta}(x_0,\nabla \varphi(x_0),-D_\gamma u(x_0))].
    \end{array}
    \end{equation}
where
$G_{x_0}^{\alpha\beta}=G^{\alpha\beta}(D_\gamma(u-\varphi)(x_0)\mathfrak{C}_{\alpha\beta}(x_0)+\nabla_{\alpha\beta}\varphi(x_0)
            -\mathcal{A}_{\alpha\beta}(x_0,\nabla \varphi(x_0),-D_\gamma u(x_0)))$.

We consider two possible cases:

{\it Case 1.} $G[u]\geq G[\underline u]/2$ at $x_0$. An upper bound
for $D_{nn}u(x_0)$ follows directly from equation \eqref{3.25}.

{\it Case 2.} $G[u]<G[\underline u]/2$ at $x_0$. Fixing a principle
coordinate system at the point $x_0$ and a corresponding
neighbourhood $\mathcal{N}$ of $x_0$ with $\gamma_n<0$ on
$T=\mathcal{N}\cap
\partial \Omega$. Since the orthogonal vector field $\xi^{(1)},\cdots,\xi^{(n-1)}$
agrees with the coordinate system at the point $x_0$, we have
$\xi^{j}_i(x_0)=\delta_{ij}$, $i,j=1,\cdots,n-1$, $-D_\gamma
u(x_0)=D_nu(x_0)$,
$\mathfrak{C}_{\alpha\beta}(x_0)=D_\alpha\gamma_\beta(x_0)$, and
$\mathcal{A}_{\alpha\beta}(x_0,\nabla \varphi(x_0),-D_\gamma
u(x_0))=A_{\alpha\beta}(x_0,D'\varphi(x_0),D_nu(x_0))$, where
$D'=(D_1,\cdots,D_{n-1})$. By the concavity of $G$, we have, at
$x_0$,
    \begin{equation}\label{3.37}
        \begin{array}{ll}
            & G[u](x_0)-G[\underline u](x_0)\\
            \geq & G_{x_0}^{\alpha\beta}\left[ D_\gamma(u-\varphi)(x_0)\mathfrak{C}_{\alpha\beta}(x_0)+\nabla_{\alpha\beta}\varphi(x_0)
            -\mathcal{A}_{\alpha\beta}(x_0,\nabla \varphi(x_0),-D_\gamma
            u(x_0))\right] \\
            &-G_{x_0}^{\alpha\beta}\left[D_\gamma(\underline u-\varphi)(x_0)\mathfrak{C}_{\alpha\beta}(x_0)+\nabla_{\alpha\beta}\varphi(x_0)
            -\mathcal{A}_{\alpha\beta}(x_0,\nabla \varphi(x_0),-D_\gamma
            \underline u(x_0))\right]\\
            = & -G_{x_0}^{\alpha\beta}\left[D_n(u-\underline
            u)(x_0)D_\alpha\gamma_\beta(x_0)+A_{\alpha\beta}(x_0,D'
            \varphi(x_0),D_n
            u(x_0))-A_{\alpha\beta}(x_0,D' \varphi(x_0),D_n
            \underline u(x_0))\right]\\
            \geq & -D_n(u-\underline u)(x_0)G_{x_0}^{\alpha\beta}\left[D_\alpha\gamma_\beta(x_0)+D_{p_n}A_{\alpha\beta}(x_0,D'
            \varphi(x_0),D_n
            u(x_0))\right].
        \end{array}
    \end{equation}
here the regular condition \eqref{A3W} of $A$ is used. By the
ellipticity of the strict subsolution $\underline u$, we can fix a
positive constant $\delta_0$ for which
    \begin{equation}\label{3.36}
        G[\underline u]\geq \delta_0,\hspace{2mm}{\rm for}\hspace{1.5mm}{\rm
        all}\hspace{1.5mm}x\in T.
    \end{equation}
We also have
    \begin{equation}\label{3.38}
        0< D_n(u-\underline u)(x_0)\leq \kappa,
    \end{equation}
for a positive constant $\kappa$. Combining (\ref{3.37}),
(\ref{3.36}) and (\ref{3.38}), since $G[u]<G[\underline u]/2$ at
$x_0$, we have
    \begin{equation}
        G_{x_0}^{\alpha\beta}\left[D_\alpha\gamma_\beta(x_0)+D_{p_n}A_{\alpha\beta}(x_0,D'\varphi(x_0),D_n
        u(x_0))\right]\geq \frac{\delta_0}{2\kappa}>0.
    \end{equation}
Let
    \begin{equation}
        \vartheta(x):=G_{x_0}^{\alpha\beta}\left[\mathfrak{C}_{\alpha\beta}(x)+D_{p_n}\mathcal{A}_{\alpha\beta}(x,\nabla \varphi(x),
        -D_\gamma u(x_0))\right],
    \end{equation}
we have, at $x_0$,
    \begin{equation}
        \vartheta(x_0):=G_{x_0}^{\alpha\beta}\left[D_\alpha\gamma_\beta(x_0)+D_{p_n}A_{\alpha\beta}(x_0,D'\varphi(x_0),D_n
        u(x_0))\right]\geq \frac{\delta_0}{2\kappa}>0.
    \end{equation}
Since $\vartheta(x)$ is smooth near $\partial \Omega$, we can have
    \begin{equation}\label{coefficient>0}
        \vartheta(x)\geq c>0,\ \ {\rm on} \ T,
    \end{equation}
holds for some small positive constant $c$. By the regular condition \ref{A3W}
of $A$, we observe that $\mathcal{A}_{\alpha\beta}$ is convex with respect to
$p_n$. Therefore, we have
    \begin{equation}\label{Aconvex}
        \begin{array}{ll}
            &\mathcal{A}_{\alpha\beta}(x,\nabla \varphi(x), -D_\gamma
            u(x_0))-\mathcal{A}_{\alpha\beta}(x,\nabla \varphi(x), -D_\gamma
            u(x))\\
            \leq & D_{p_n}\mathcal{A}_{\alpha\beta}(x,\nabla \varphi(x), -D_\gamma
            u(x_0))(D_\gamma u(x)-D_\gamma
            u(x_0)).
        \end{array}
    \end{equation}
Since $G_{x_0}^{\alpha\beta}$ is positive definite and
$\gamma_n(x)<0$, we have, by \eqref{Gconcavity} and \eqref{Aconvex},
    \begin{equation}\label{Dnu}
        \begin{array}{ll}
            &G_{x_0}^{\alpha\beta}\left[\mathfrak{C}_{\alpha\beta}(x)+D_{p_n}\mathcal{A}_{\alpha\beta}(x,\nabla \varphi(x),
            -D_\gamma u(x_0))\right]D_nu(x)\\
            \leq &
            \gamma_n(x)\left\{G_{x_0}^{\alpha\beta}\left[D_\gamma \varphi(x)\mathfrak{C}_{\alpha\beta}(x)+D_{p_n}\mathcal{A}_{\alpha\beta}(x,\nabla \varphi(x),
            -D_\gamma u(x_0))D_\gamma u(x_0)+\nabla_{\alpha\beta}\varphi(x_0)-\nabla_{\alpha\beta}\varphi(x)\right.\right.\\
            & \left.\left.+D_\gamma(u-\varphi)(x_0)\mathfrak{C}_{\alpha\beta}(x_0)+\mathcal{A}_{\alpha\beta}(x,\nabla \varphi(x),-D_\gamma
            u(x_0))-\mathcal{A}_{\alpha\beta}(x_0,\nabla \varphi(x_0),-D_\gamma u(x_0))\right]\right\}.
        \end{array}
    \end{equation}
Observing that $\vartheta(x)$ is the coefficient of $D_nu(x)$ in
\eqref{Dnu}, by \eqref{coefficient>0}, we have
    \begin{equation}\label{Dnu<Theta}
        D_nu(x)\leq \Theta(x),\ \ {\rm on} \ T,
    \end{equation}
where $\Theta(x)$ is a smooth function defined by
    \begin{equation}
        \begin{array}{c}
            \Theta(x):=\gamma_n(x)(\vartheta(x))^{-1}\left\{G_{x_0}^{\alpha\beta}\left[D_\gamma \varphi(x)\mathfrak{C}_{\alpha\beta}(x)+D_{p_n}\mathcal{A}_{\alpha\beta}(x,\nabla \varphi(x),
            -D_\gamma u(x_0))D_\gamma u(x_0)+\nabla_{\alpha\beta}\varphi(x_0)\right.\right.\\
            \left.\left.-\nabla_{\alpha\beta}\varphi(x)+D_\gamma(u-\varphi)(x_0)\mathfrak{C}_{\alpha\beta}(x_0)+\mathcal{A}_{\alpha\beta}(x,\nabla \varphi(x),-D_\gamma
            u(x_0))-\mathcal{A}_{\alpha\beta}(x_0,\nabla \varphi(x_0),-D_\gamma u(x_0))\right]\right\}.
        \end{array}
    \end{equation}
We now define a function on $T$,
    \begin{equation}
        w(x):=D_nu(x)- \Theta(x).
    \end{equation}
By extending $\varphi$, $\gamma$ and $\mathfrak{C}$ smoothly to the
interior near the boundary to be constant in the normal direction,
the function $w$ is extended to $\Omega_\delta=\Omega\cap
B_\delta(x_0)$ in the neighbourhood of $x_0$. By calculation, we
have
    \begin{equation}
        |\mathcal{L}w|=|\mathcal{L}D_nu(x)-\mathcal{L}\Theta(x)|\leq C(1+\sum_i
        F^{ii}),\ \ {\rm in} \ \Omega_\delta,
    \end{equation}
here the differentiated equation \eqref{LDku} for $k=n$ is used.
Since $w$ is extended to be constant in the normal direction, we
have, by \eqref{Dnu<Theta},
    \begin{equation}
        w\leq 0,\ \ {\rm on} \ \partial \Omega_\delta.
    \end{equation}
Therefore, we obtain
    \begin{equation}
    \left\{\begin{array}{rllll}
            \mathcal{L}(K'\psi \pm w) & \leq & 0, & {\rm in} &
            \Omega_\delta,\\
            K'\psi-w & \geq & 0, & {\rm on} & \partial \Omega_\delta,
        \end{array}\right.
    \end{equation}
where $K'$ is a sufficiently large constant, $\psi$ is the barrier
function constructed in Lemma \ref{Lemma 3.1}. By the maximum
principle, we have
    \begin{equation}
        K'\psi-w \geq 0, \ \  {\rm in} \  \Omega_\delta,
    \end{equation}
which leads to
    \begin{equation}
        D_n(K'\psi-w)\geq 0,\ \ {\rm at}\ x_0,
    \end{equation}
namely
    \begin{equation}
        D_{nn}u(x_0)\leq C.
    \end{equation}
Utilizing equation \eqref{3.25} again, we conclude from the above
two cases a positive lower bound from below for $G[u]$ at $x_0$.
Recall that $G[u]$ attains its minimum at $x_0$. Hence, by
\eqref{3.25} we finally obtain an upper bound
$D_{\gamma\gamma}u(x)\leq C$ at any boundary point $x\in \partial
\Omega$.

By the ellipticity, we have $tr(D^2u-A(x,Du))\geq 0$, that is
    \begin{equation}
        D_{nn}u\geq \sum_i^n A_{ii}(x,Du)-\sum_i^{n-1}D_{ii}u,
    \end{equation}
we now get an estimate from below for $D_{nn}u$ on $\partial
\Omega$.

In conclusion, we have obtained the following desired second
derivative bound on the boundary,
    \begin{equation}\label{boundary bound}
        \sup_{\partial \Omega} |D^2u|\leq C,
    \end{equation}
where the constant $C$ depends on $n, A, B, \Omega, \varphi,
\underline u$ and $\sup\limits_\Omega{(|u|+|Du|)}$.

We see that the second derivative estimate \eqref{boundary bound} on
the boundary holds for solutions of augmented Hessian equations
\eqref{1.1} for all $1 \leq k \leq n$. In fact, if we only consider
the semilinear case when $k=1$, the proof will be much simpler. In
this case, the double tangential derivative bound can also be
derived from \eqref{Dalphabeta}, and the double normal derivative
bound can be obtained directly from the equation $\Delta
u-\sum\limits_i^n A_{ii}(x,Du)=B(x,Du)$.

Now, we can easily give the proof of Theorem 1.2.

\vspace{3mm}

\noindent \textbf{Proof of Theorem 1.2}. By the global estimate
\eqref{interior bound} and the above boundary estimate
\eqref{boundary bound}, we obtain the {\it a priori} second
derivative estimate \eqref{second bound} and Theorem \ref{Th2}
follows.

\qed

\begin{remark}

In order to obtain the second derivative estimates for admissible
solutions of \eqref{1.1} on the boundary without the subsolution
assumption, a kind of geometric condition should be imposed on the
domain $\Omega$. Recall that in the case for Monge-Amp\`ere type
equation \cite{JTY2013,Tru2006}, the concept of domain $A$-convexity
is introduced extending that of $c$-convexity in optimal
transportation. If $\Omega$ is a connected domain in $\mathbb{R}^n$
with $\partial \Omega \in C^2$, and $A\in
C^1(\Omega\times\mathbb{R}^n;\mathbb{S}^n)$, we say that $\Omega$ is
uniformly $(k-1)$-$A$-convex with respect to $u$, if
    \begin{equation}\label{Akconvex}
        S_{k-1}(\kappa)\geq \delta_0 > 0,
    \end{equation}
for some $\delta_0>0$ and all $p\in Du(\Omega)$,
where $\kappa=(\kappa_1,\cdots,\kappa_{n-1})$
denote the eigenvalues of
    \begin{equation}\label{eigenvalues}
        \{D_i\gamma_j(x)-A_{ij,p_k}(x,p)\gamma_k(x)\},
    \end{equation}
for all $x\in \partial \Omega$, and unit outer normal $\gamma$. In
other words, a domain $\Omega$ is called uniformly
$(k-1)$-$A$-convex if the eigenvalues of the matrix
$\{D_i\gamma_j(x)-A_{ij,p_k}(x,p)\gamma_k(x)-c_0\delta_{ij}\}$ (as a
vector in $\mathbb{R}^{n-1}$) lie in $\Gamma_{k-1}$, provided $c_0$
is a sufficiently small positive constant, where $\delta_{ij}$ is
the usual Kronecker delta. When $k=n$, the $(n-1)$-$A$-convexity
corresponds to the $A$-convexity in \cite{JTY2013,Tru2006}. When
$A\equiv 0$, the $(k-1)$-$A$-convexity corresponds to the
$(k-1)$-convexity as in \cite{Tru1997,Wang2009}. Particularly, when
$k=n$ and $A\equiv 0$, the $(k-1)$-$A$-convexity reduces to the
usual convexity. If we impose a uniformly $(k-1)$-$A$-convexity
condition on the domain $\Omega$, we can still obtain the second
derivative estimates on the boundary for admissible solutions of
\eqref{1.1}. We remark that the subsolution assumption can be
replaced by the $A$-boundedness condition and the uniformly
$(k-1)$-$A$-convexity condition of the domain $\Omega$ in the main
theorems of this paper, but in some sense it is a bit weaker to
assume the existence of a subsolution with the same boundary trace.
\end{remark}


\section{ General type equations}\label{Section 5}

In many applications, the matrix function $A$ and the right hand
side term $B$ in the augmented Hessian equations may also depend on $u$.
In this section, we consider some more general type augmented Hessian equations with the form
\eqref{general} and formulate the corresponding theorems for
the {\it a priori} $C^2$ estimates of the solutions. The general type equation \eqref{general} reduces to the
Monge-Amp\`ere type equation when $k=n$, which has applications in near field
optics, see \cite{KW2010,KO1997,LiuTru2013,Tru2012}.

With both the matrix function $A$ and the scalar
function $B$ depending on $u$, we still use similar notation as in
the previous sections. Here the augmented Hessian matrix is
$\{w_{ij}\}=\{u_{ij}-A_{ij}(x,u,Du)\}$ and the right hand side term
is $B(x,u,Du)$. We rewrite the equation \eqref{general} in the following form
    \begin{equation}\label{re general}
        F[u]=(S_k)^{\frac{1}{k}}[D^2u-A(x,u,Du)]=\tilde B(x,u,Du), \ \ {\rm in}\ \Omega,
    \end{equation}
where $\tilde B = B^{\frac{1}{k}} >0$ in $\bar \Omega \times \mathbb{R} \times \mathbb{R}^n$. The corresponding Dirichlet boundary condition is
    \begin{equation}\label{re boundary}
        u=\varphi,\ \ {\rm on} \ \partial \Omega,
    \end{equation}
where $\varphi$ is a smooth function on $\partial \Omega$. The matrix function $A$ is regular if
    \begin{equation}\label{A3W general}
        A_{ij,kl}(x,u,p)\xi_i\xi_j\eta_k\eta_l\geq 0,
    \end{equation}
for all $(x,u,p)\in\bar\Omega\times\mathbb{R}\times\mathbb{R}^n$,
$\xi,\eta\in\mathbb{R}^n$, $\xi\perp\eta$, where
$A_{ij,kl}=D^2_{p_kp_l}A_{ij}$.

In order to construct the global barrier function as in Lemma
\ref{Lemma 2.1}, we need to assume additional monotonicity conditions on both $A$
and $\tilde B$ with respect to $u$. We assume the matrix function $A(x,u,p)$ is monotone with respect to $u$, that is
    \begin{equation}\label{monotone}
        D_uA_{ij}(x,u,p)\xi_i\xi_j\geq 0,
    \end{equation}
for all $(x,u,p)\in \Omega\times\mathbb{R}\times\mathbb{R}^n$,
$\xi\in \mathbb{R}^n$. We also assume the scalar function $\tilde B(x,u,p)$ is monotone with respect to $u$, that is
    \begin{equation}\label{B_u>0}
        \tilde B_u(x,u,p)\geq 0,
    \end{equation}
for all $(x,u,p)\in \Omega\times\mathbb{R}\times\mathbb{R}^n$.
Note that these monotonicity conditions can guarantee the comparison principle and the uniqueness of the solution.
To obtain the global second derivative bound, we assume that $\tilde B(x,u,p)$ is convex with respect to $p$, that is
    \begin{equation}\label{convexity on B'}
        D^2_{p_kp_l}\tilde B(x,u,p)\xi_k\xi_l\geq 0,
    \end{equation}
for all $(x,u,p)\in \Omega\times\mathbb{R}\times\mathbb{R}^n$, $\xi\in \mathbb{R}^n$.

We define the linearized operators of $F$ by
    \begin{equation}
        L=F^{ij}[D_{ij}-D_{p_k}A_{ij}(x,u,Du)D_k],
    \end{equation}
and
    \begin{equation}
        \mathcal {L}=L-\tilde B_{p_i}D_i=F^{ij}[D_{ij}-D_{p_k}A_{ij}(x,u,Du)D_k]-\tilde B_{p_i}D_i.
    \end{equation}
where $F^{ij}=\frac{\partial F}{\partial w_{ij}}$. With the monotonicity conditions \eqref{monotone} and \eqref{B_u>0}, we can construct the barrier function similar to Lemma \ref{Lemma 2.1} by assuming the existence of an admissible supersolution.

\begin{lemma}\label{Lemma 6.2}
Assume the function $A(x,u,p)\in
C^2(\bar\Omega\times\mathbb{R}\times\mathbb{R}^n)$ is regular, $\tilde B(x,u,p)\in
C^2(\bar\Omega\times\mathbb{R}\times\mathbb{R}^n)$,
$A$ and $\tilde B$ satisfy the monotonicity conditions \eqref{monotone} and \eqref{B_u>0} respectively.
Suppose $u\in C^2(\bar \Omega)$ is an elliptic solution of equation
\eqref{re general} associated with Dirichlet the boundary condition \eqref{re boundary}, $\bar u\in C^2(\bar \Omega)$ is an admissible strict
supersolution of equation \eqref{re general} satisfying
    \begin{equation}
        (S_k)^{\frac{1}{k}}[D^2 \bar u - A(x,\bar u, D\bar u)] \leq \tilde B(x,\bar u, D\bar u)-\delta, \ \ {in} \ \Omega,
    \end{equation}
for some positive constant $\delta$, and $\bar u \geq \varphi$ on $\partial \Omega$. If $A$ is regular satisfying \eqref{A3W general}, then
    \begin{equation}
        \mathcal {L} \left({e^{K(\bar u-u)}}\right) \geq
        \epsilon_1\sum\limits_{i} F^{ii}-C,
    \end{equation}
holds in $\Omega$ for positive constants $K$, $\epsilon_1$ and $C$, which depend on $A, \tilde B, \Omega, \|u\|_{C^1}$
and $\|\bar u\|_{C^1}$.
\end{lemma}
\noindent \textbf{Proof}. Since $\bar u$ is an admissible strict supersolution
of \eqref{re general}, for any $x_0\in \Omega$, the perturbation
function $\bar u_\epsilon=\bar u -\frac{\epsilon}{2}|x-x_0|^2$ is
still an admissible strict supersolution satisfying
    \begin{equation}
        F[\bar u_\epsilon]=(S_k)^{\frac{1}{k}}[D^2\bar u_\epsilon-
        A(x,\bar u_\epsilon,D\bar u_\epsilon)]\leq \tilde B(x,\bar u_\epsilon,D\bar u_\epsilon)-\tau, \ \ {\rm in}\ \Omega,
    \end{equation}
for some positive constant $\tau$. Consequently, we have
    \begin{equation}\label{0 leq}
        0 \leq F[\bar u_\epsilon] < \tilde B(x,\bar u_\epsilon,D\bar u_\epsilon), \ \ {\rm in}\ \Omega.
    \end{equation}
By the monotone assumptions \eqref{monotone} and \eqref{B_u>0}, we can automatically have $u < \bar u_\epsilon$ in $\Omega$ from the comparison principle.

Let $v=\bar u-u$, $v_\epsilon=\bar
u_\epsilon-u$. By calculation, we get
    \begin{equation}\label{6.6}
        \begin{array}{lll}
            Lv&=& L(v_\epsilon)+L(\frac{\epsilon}{2}|x-x_0|^2)\\
            &=& \epsilon F^{ii}-\epsilon
            F^{ij}D_{p_k}A_{ij}(x,u,Du)(x-x_0)_k\\
            & & +F^{ij}\{D_{ij}(\bar u_\epsilon -u)-[A_{ij}(x,\bar
            u_\epsilon,D\bar
            u_\epsilon)-A_{ij}(x,u,Du)]\}\\
            & & +F^{ij}\{A_{ij}(x,u,D\bar u_\epsilon)-A_{ij}(x,u,Du)-D_{p_k}A_{ij}(x,u,Du)D_kv_\epsilon\}\\
            & & +F^{ij}\{A_{ij}(x,\bar u_\epsilon,D\bar
            u_\epsilon)-A_{ij}(x,u,D\bar u_\epsilon)\}.
        \end{array}
    \end{equation}
Since $u \leq \bar u_\epsilon$ and $A$ satisfies the monotonicity condition \eqref{monotone}, we have
    \begin{equation}
        F^{ij}\{A_{ij}(x,\bar u_\epsilon,D\bar
            u_\epsilon)-A_{ij}(x,u,D\bar
            u_\epsilon)\}=F^{ij}D_uA_{ij}(x,\hat u, D\bar
            u_\epsilon)(\bar u_\epsilon-u)\geq 0,
    \end{equation}
where $\hat u=\theta \bar u_\epsilon + (1-\theta)u$ for some
$\theta\in (0,1)$.

By the concavity of $F$, we have
    \begin{equation}
        F^{ij}\{D_{ij}(\bar u_\epsilon -u)-[A_{ij}(x,\bar
            u_\epsilon,D\bar
            u_\epsilon)-A_{ij}(x,u,Du)]\}\geq F[\bar
            u_\epsilon]-F[u]= F[\bar
            u_\epsilon]- \tilde B(x,u,Du)\geq - C_1,
    \end{equation}
where the first inequality of \eqref{0 leq} is used, and $C_1$ is a positive constant depending on $\tilde B$ and
$\sup\limits_\Omega{(|u|+|Du|)}$.

By the regular condition \eqref{A3W general}, other terms in \eqref{6.6} can be estimated in the same manner as in the proof of Lemma
\ref{Lemma 2.1}. We only need to substitute $\bar u$ and $\bar u_\epsilon$ in place of $\underline
u$ and $\underline u_\epsilon$ respectively in that proof.

Then by the finite covering, the conclusion of Lemma \ref{Lemma 6.2} can be completed by following the steps in the
proof of Lemma \ref{Lemma 2.1}.

\qed

Similar to Theorem \ref{Th1} and Theorem \ref{Th2}, we can use the barrier
function in Lemma \ref{Lemma 6.2} to derive the second order
derivative estimates for solutions of \eqref{re general}-\eqref{re boundary}.
We formulate these results as follows:
\begin{theorem}\label{Th5.1}
In addition to the assumptions in Lemma \ref{Lemma 6.2}, suppose the condition \eqref{convexity on B'} holds for $\tilde B$. Then we have the estimate
    \begin{equation}
        \sup\limits_\Omega |D^2u| \leq C(1+\sup\limits_{\partial
        \Omega} |D^2u|),
    \end{equation}
where the constant $C$ depends on $n, A, \tilde B, \Omega, \bar u$
and $\sup\limits_\Omega{(|u|+|Du|)}$.
\end{theorem}

\begin{remark}\label{remark 5.1} In the paper \cite {JT2014} we obtain corresponding second derivative bounds for the special case of generated prescribed Jacobian equations under the conditions G1,G2, G1*, G3w, G4w introduced in \cite {Tru2012}, (without convexity or monotonicity hypotheses on $B$).
These extend the optimal transportation case as indicated  in Remark 2.4 and depend on the construction of a suitable barrier $\underline u$.
\end{remark}

In order to obtain the second order derivative estimates on the boundary, we still need a subsolution $\underline u$ with the boundary value $\varphi$ as the deduction in Section \ref{Section 4}. The subsolution is necessary to derive the double normal derivative bound on the boundary. Combining the boundary second derivative estimates with Theorem \ref{Th5.1}, we have the following result for the second
order {\it a priori} estimates for Dirichlet problem \eqref{re general}-\eqref{re boundary} of the general type equations.
\begin{theorem}\label{Th5.2}
In addition to the assumptions in Theorem \ref{Th5.1}, suppose also there exists an elliptic
subsolution $\underline u\in C^2(\bar\Omega)$ of the equation \eqref{re general} with $\underline u = \varphi$ on $\partial \Omega$. Then we
have the global estimate,
    \begin{equation}
        \sup_\Omega |D^2u| \leq C,
    \end{equation}
where the constant $C$ depends on $n, A, \tilde B, \Omega, \varphi, \bar u$, $\underline u$ and $\sup\limits_\Omega{(|u|+|Du|)}$.
\end{theorem}

We have obtained the global second derivative estimates for solutions of the Dirichlet problem \eqref{re general}-\eqref{re boundary}.
If the lower order derivative estimates are also obtained,
we can obtain the $C^{2,\alpha}$ estimates of the solutions by the Evans-Krylov theorems,\cite{GT2001}.
Using the method of continuity, we can have the existence and uniqueness of the
classical solutions of the Dirichlet problem \eqref{re general}-\eqref{re boundary}.

\begin{remark}\label{remark 6.2}
We remark that the monotonicity condition \eqref{B_u>0}
for right hand side function $\tilde B$ is to guarantee the
comparison principle and the uniqueness of the solution. Since we
use the subsolution $\underline u$ and the supersolution $\bar u$
here, we need to compare them with the solution $u$. If we assume
the geometric conditions \eqref{Abound} and \eqref{Akconvex} of the
domain $\Omega$ rather than the existence of a subsolution,
condition \eqref{B_u>0} is not necessary in Theorem \ref{Th5.2} and
may only be used to guarantee the uniqueness by the continuity method.
Furthermore, if the condition $\tilde B_u\geq 0$ is not
satisfied, one may still have the existence theorem for the
classical solutions by using Leray-Schauder fixed point theorem (see
Theorem 11.6 in \cite{GT2001}) or a degree argument used in
\cite{CNSIII,Guanbo2007}. We refer the reader to the proof of
Theorem 1.5 in \cite{Urbas1995} for the procedure of using
the Leray-Schauder fixed point theorem for similar purpose.
\end{remark}


\end{document}